\documentstyle{amsppt}
\magnification=\magstep1
\pageheight{23truecm}
\NoBlackBoxes

\def\Z{{\Bbb Z}}
\def\R{{\Bbb R}}
\def\Sb{{\Bbb S}}
\def\T{{\Bbb T}}
\def\B{{\Bbb B}}
\def\noi{\noindent}
\def\lk{\operatorname{lk}}
\def\cl{\operatorname{cl}}
\def\rank{\operatorname{rank}}
\def\ss{\smallskip}
\def\ms{\medskip}
\def\bs{\bigskip}
\def\kn{\mathaccent"7017 }
\def\kw{\mathbin{\raise 1pt \hbox {$\scriptstyle\rlap
{$\scriptstyle \sqcap $} \sqcup $}}}
\catcode`@=11
\def\dwiestrz#1#2{\,\vcenter{\m@th\ialign{##\crcr
      $\hfil\scriptstyle{\ #1\ }\hfil$\crcr
      \noalign{\kern2\p@\nointerlineskip}
      \rightarrowfill\crcr
      \noalign{\kern\p@\nointerlineskip}
      \leftarrowfill\crcr\noalign{\kern2\p@\nointerlineskip}
      $\hfil\scriptstyle{\ #2\ }\hfil$\crcr}}\,}
\catcode`@=13

\rightheadtext{Generalized manifolds in products of curves}
\leftheadtext{A. Koyama, J. Krasinkiewicz and S. Spie{\D z}}

\topmatter


\null\vskip1cm
\title
Generalized manifolds in products of curves \vskip-3mm
\endtitle

\author
Akira Koyama, J{\' o}zef Krasinkiewicz, and Stanis{\l}aw Spie{\D
z} \vskip2mm
\endauthor

\address
Department of Mathematics, Faculty of Science, Suruga, Shizuoka
University, Shizuoka, 422-8529, Japan
\endaddress
\email sakoyam$\@$ipc.shizuoka.ac.jp
\endemail

\address
The Institute of Mathematics, Polish Academy of Sciences, ul.{\'
S}niadeckich 8, 00-950, Warsaw, Poland \endaddress
\address
Institute of Mathematics and Informatics, University of Opole, ul.
Oleska 48, 45-052 Opole, Poland
\endaddress
\email jokra$\@$impan.gov.pl
\endemail

\address
The Institute of Mathematics, Polish Academy of Sciences, ul.{\'
S}niadeckich 8, 00-950, Warsaw, Poland
\endaddress
\email spiez$\@$impan.gov.pl
\endemail

\abstract The intent of this article is to study some special
$n$-dimensional continua lying in products of $n$ curves. (The
paper may be viewed as a sequel to \cite{K-K-S 1}; majority of its
content is an improved version of a portion of \cite{K-K-S}.) We
show that if $X$ is a locally connected, so-called, quasi
$n$-manifold lying in a product of $n$ curves then $\rank
H^1(X)\ge n$. Moreover, if $\rank H^1(X)<2n$ then $X$ can be
represented as a product of an $m$-torus and a quasi
$(n-m)$-manifold, where $m\ge2n-\rank H^1(X)$. It follows that
certain 2-dimensional contractible polyhedra are not embeddable in
products of two curves. On the other hand, we show that any
collapsible 2-dimensional polyhedron can be embedded in a product
of two trees. We answer a question of R. Cauty proving that closed
surfaces embeddable in products of two curves can be also embedded
in products of two graphs. On the other hand, we construct an
example of a 2-dimensional polyhedron which can be embedded in a
product of two curves though it is not embeddable in any product
of two graphes. This solves in the negative another problem of
Cauty. Some open problems have been included.
\endabstract



\endtopmatter

\vfill \break

\null \vskip2cm \bs\centerline{\bf CONTENTS}

\bs \centerline{\bf 1. Introduction}

\bs \centerline{\bf 2. Generalized manifolds}

\bs \item{\eightrm 2A.} {\eightit Definitions and basic properties
of certain generalized manifolds \hfill{6}}

\item{\eightrm 2B.} {\eightit On locally connected quasi
manifolds. From embeddings into products of curves \hfill
\vskip-1mm to embeddings into products of graphs\hfill{10}}

\item{\eightrm 2C.} {\eightit  Ramified manifolds in products of
graphs\hfill{14}}

\item{\eightrm 2D.} {\eightit Product structure of generalized
manifolds lying in products of curves\hfill{17}}

\item{\eightrm 2E.} {\eightit Contractible {\eightrm
2}-dimensional polyhedra in products of two graphs\hfill{20}}

\item{\eightrm 2F.} {\eightit On embeddings of {\eightrm
2}-dimensional polyhedra into products of two graphes \hfill
\vskip-1mm -- a solution of Cauty's problem\hfill{23}}

\bs \centerline {\bf 3. Problems}

\vfill \break

\document

\head 1. Introduction
\endhead

Throughout this paper we use the following standard conventions,
terminology and notation. All {\it spaces} under discussion are
assumed to be metrizable and all {\it mappings} (also called {\it
maps}) -- continuous. By a {\it compactum} we mean a compact
(metric) space, by a {\it continuum} we mean a connected
compactum, and by a {\it curve} we mean a 1-dimensional continuum.
A continuum $X$ is said to be {\it unicoherent} if for every
presentation of $X$ as a union of two continua the intersection of
these continua is connected. We write $X\approx Y$ to indicate
that $X$ is homeomorphic to $Y$.

By $\B^n$, $n\geq1$, we denote the closed unit $n$-ball in the
Euclidean $n$-space $\R^n$. A space homeomorphic to $\B^n$ is
called a (closed) $n$-{\it disc}. Sometimes $2$-discs will be
simply called {\it discs}. By $\Sb^{n-1}$ we denote the unit
$n$-sphere in $\R^n$ - the boundary of $\B^n$. A space
homeomorphic to $\Sb^n$ is called a {\it topological $n$-sphere};
a space homeomorphic to $\Sb^1$ is called a {\it $($topological$)$
circle} (or a {\it simple closed curve}). A space homeomorphic to
$\Sb^1\times I$, where $I=[0, 1]$ stands for the unit interval in
$\R$, is called a {\it $($topological$)$ cylinder}. A space
homeomorphic to the open unit $n$-ball $\kn\B^n=\B^n \backslash
\Sb^{n-1}$ is called an {\it open $n$-disc}. As usual, by the
$n$-torus $\T^n$, $n\geq1$, we mean the $n$-fold product
$\Sb^1\times\cdots\times \Sb^1$. In particular, $\T^1=\Sb^1$. By
$\T^0$ we mean a one-point space. A space homeomorphic to $\T^n$
is called a {\it topological n-torus}.

By a {\it polyhedron} we mean the underlying space $|K|$ of a
finite regular $CW$ complex $K$. Recall that a $CW$ complex $K$ is
said to be {\it regular} if each cell $\sigma\in K$ admits a
characteristic map $\varphi_{\sigma}: {\B}^n \to \sigma$, where
$n=\dim \sigma$, which is a homeomorphism. Any {\it simplicial
complex} $K$ is regular, such $K$ is called a {\it triangulation}
of $|K|$. By a {\it graph} we mean a 1-dimensional polyhedron, and
by a {\it tree} we mean a connected graph containing no simple
closed curve. By a {\it dendrite} we mean a non-degenerate locally
connected continuum containing no simple closed curve. A
non-degenerate continuum is said to be a {\it local dendrite}
({\it local tree}, respectively) if every point has a closed
neighborhood which is a dendrite (tree, respectively). It is known
that dendrites coincide with 1-dimensional compact absolute
retracts, and local dendrites -- with 1-dimensional compact
absolute neighborhood retracts (cf. \cite{Kur}).

In 1958 J. Nagata \cite{N1} discovered the following remarkable
theorem.

\proclaim{Theorem 1.1 (Nagata)} Every $n$-dimensional space,
$n\ge2$, can be embedded in the topological product
$X_1\times\cdots \times X_{n+1}$ of $1$-dimensional spaces.\qed
\endproclaim

Nagata also asked whether $n$ 1-dimensional spaces would suffice
\cite{N2, p.163}. This question and the theorem itself gave rise
to interesting investigations. The question itself was answered in
the negative by K.~Borsuk in $1975$  \cite{Bo3}. Actually, Borsuk
proved the following interesting result.

\proclaim{Theorem 1.2 (Borsuk)} The $n$-sphere $\Sb^n$, $n\ge 2$,
is not embeddable in any product of n curves.\qed
\endproclaim

In our paper we shall see that some 2-dimensional contractible
(so, acyclic) polyhedra have this property as well.

An $n$-dimensional space, $n\geq2$, is said to be {\it ordinary}
if it can be embedded in a product of $n$ 1-dimensional spaces;
otherwise we call it {\it exceptional}. Hence $\Sb^2$ and some
other spaces are exceptional. Any 1-dimensional compactum can be
embedded in the Menger curve $\mu$. It follows that an
$n$-dimensional compactum is ordinary if and only if it can be
embedded in the $n$-fold product $\mu^n$.

All {\it manifolds} discussed in this paper are assumed to be
compact and connected (possibly with non-empty boundary), unless
opposite is explicitly stated. 2-manifolds are often called {\it
surfaces}. The interior of $M$ will be often denoted by $\kn M$. A
manifold $M$ is {\it closed} if its boundary is empty,
$\partial{}M=\emptyset$. A mapping $f:X\to{M}$, where $M$ is a
closed manifold, is said to be {\it essential} if every mapping
$g:X\to{M}$ homotopic to $f$ is surjective.

The symbol $H^{\ast}(\cdot\,;\,G)$ is used to denote the {\v Cech}
cohomology functor with coefficients in an Abelian group $G$. In
some cases where no confusion is likely to occur we shall write
$f^{\ast}$ instead of $H^{\ast}(f\,;\,G)$, where $f$ is a mapping.
The cohomology functor $H^{\ast}(\cdot\,,\,\Z)$ with integer
coefficients $\Z$ will be abbreviated to $H^{\ast}(\cdot)$. Thus
the groups $H^{\ast}(X\,,\,\Z)$ and the homomorphisms
$H^{\ast}(f\,,\,\Z)$ will be written briefly $H^{\ast}(X)$ and
$H^{\ast}(f)$, respectively. By the Hopf Classification Theorem,
cf. \cite{Sp, p.431}, for any $n$-dimensional space $X$, the group
$H^n (X)$ is in one-to-one correspondence with the set of homotopy
classes of maps $X\to \Sb^n$. Non-zero elements correspond to
homotopy classes of essential maps. As usual, $H_*(X)$ denotes the
singular homology functor with integer coefficients.

Let $g_1,\cdots,g_k$ be elements of an Abelian group $G$. They are
said to be {\it linearly independent} (over $\Z$) if the equality
$n_1g_1+\cdots+n_kg_k=0$, $n_i\in\Z$, implies $n_1=\cdots=n_k=0$.
By the {\it rank} of $G$, denoted $\rank G$, we mean the maximal
number of linearly independent elements in $G$ (over $\Z$). We
write $G\cong H$ to denote that the groups are isomorphic.

A point $x\in X$ is said to be {\it of order} $n$ if $n$ is the
minimal number such that $x$ admits arbitrarily small open
neighborhoods with boundaries containing at most $n$ points.
Points of order 1 are called {\it endpoints} of the space.

Now we present a brief summary of the main results of our paper.
An $n$-dimensional continuum $X$, $n\geq 1$, is said to be a {\it
quasi n-manifold} if for every point $x\in X$ there is an open
neighborhood $V$ of $x$ such that every closed subset of $X$ which
separates $X$ between $x$ and $X\setminus V$ and has dimension
$\leq n-1$ admits an essential map into $\Sb^{n-1}$ (see Section
2A for details). This and some other classes of $n$-dimensional
continua have been defined in Chapter 2. Each class comprises all
$n$-manifolds, and the following holds.

\ms\noi{\bf Theorem 1.3} (cf. Theorem 2B.1){\bf.} {\it Let $X$ be
a locally connected quasi $n$-manifold, $n\ge2$, with $H^1(X)$ of
finite rank. If $X$ embeds in a product of $n$ curves then there
exists an embedding $g =(g_1, \cdots ,g_n):X \to P_1 \times \cdots
\times P_n$ such that} \roster \item "(1)" {\it each $P_i$ is a
graph with no endpoint} $(${\it hence} $\rank H^1(P_i)\ge 1$$)$;
\item "(2)" {\it each $g_i$ is a monotone surjection} $(${\it
hence} $\rank H^1(P_i)\le \rank H^1(X)$$)$.\qed
\endroster

\ms\noi{\bf Corollary 1.4} (see Corollary 2B.3){\bf.} {\it If a
closed $n$-manifold is embeddable in a product of $n$ curves then
it is also embeddable in a product of $n$ graphs.} \qed

\ms\noindent It follows that {\it if a closed surface can be
embedded in a product of two curves then it can be also embedded
in a product of two graphs}. This answers a question posed by
Cauty in \cite{C}. It turns out that other analogous question
explicitly formulated in that paper - where the word "surface" is
replaced by "$2$-dimensional polyhedron" - has negative answer.
Actually, we shall construct the following example.\ms

\noi{\bf Theorem 1.5} (cf. Theorem 2F.6){\bf.} {\it There exists a
$2$-dimensional polyhedron which can be embedded in a product of
two curves though it is not embeddable in any product of two
graphs.} \qed\ms

We shall also prove a stronger version of the following

\ms\noi{\bf Theorem 1.6} (cf. Theorem 2D.6){\bf.} {\it Let $X$ be
a locally connected quasi $n$-manifold lying in a product of $n$
curves. Then $\rank H^1(X)\ge n$.} \qed\ms

A 1-cell of a $2$-dimensional regular $CW$ complex $K$ is said to
be {\it free in} $|K|$ if it is incident with exactly one $2$-cell
of $K$.

\ms\noi{\bf Corollary 1.7}{\bf.}{\it No contractible
$2$-dimensional polyhedron $|K|$ without free edges can be
embedded in a product of two curves}. \qed\ms

\noi There are two well known polyhedra satisfying the hypotheses
of this corollary: the Borsuk example \cite{Bo1} (which occurs in
\cite{Z} under the name {\it "dunce hat"}), and the {\it "Bing
house"}, cf. \cite{R-S}. Hence neither can be embedded in any
product of two curves. We shall prove also other results on
$2$-dimensional polyhedra. Here are some of them.

\ms\noi{\bf Theorem 1.8} (cf. Theorem 2E.1){\bf .} {\it Let $X$ be
a $2$-dimensional connected polyhedron. If $X$ can be embedded in
a product of two curves and $\rank{}H^1(X)\le2$, then $X$
collapses to either a point, or a graph, or a torus. In
particular, $X$ is collapsible if $\rank{}H^1(X)=0$.} \qed

\ms\noi{\bf Theorem 1.9} (cf. Theorem 2E.3) {\bf.} {\it Let $X$ be
a $2$-dimensional polyhedron. If $X$ is collapsible then $X$ can
be embedded in a product of two trees.} \qed

\ms\noi{\bf Theorem 1.10} (see Corollary 2E.7){\bf.} {\it The cone
over an $n$-dimensional polyhedron can be embedded in a product of
$n+1$ copies of an $m$-od.} \qed

\ms \noi By an $m$-{\it od} we mean the cone over an $m$-element
set.

\head 2. Generalized manifolds
\endhead

This chapter splits into sections 2A-2E. In section 2A we define
some broad classes of continua each comprising all closed
manifolds of corresponding dimension. We name them {\it
quasi-manifolds}, {\it pseudo-manifolds}, {\it para-manifolds},
and {\it ramified-manifolds} and present some basic properties of
these classes. In Theorem 2B.1 we prove that embeddings of locally
connected quasi $n$-manifolds in products of $n$ curves can be
factored through some special embeddings in products of $n$ graphs
which all are quasi 1-manifolds. Then we construct a closed
surface lying in a product of two curves whose image under either
projection is not a graph (Example 2B.2). Theorem 2B.1 has
noteworthy consequences. For example, it follows that no locally
connected and unicoherent quasi $n$-manifold can be embedded in
any product of $n$ curves. In section 2C we present a list of
basic properties of ramified $n$-manifolds lying in products of
$n$ graphs. To obtain these properties we carefully study the
"fibers" of the projections restricted to the ramified manifolds.
In section 2D we prove the fundamental Theorem 2D.6 on algebraic
and product structures of locally connected quasi $n$-manifolds
lying in products of $n$ curves. In particular, it implies that
there exist contractible 2-dimensional polyhedra not embeddable in
any product of two curves. Thus we reveal acyclic polyhedra which
share the property of $\Sb^2$ from the Borsuk Theorem 1.1. In
section 2E we prove that any $2$-dimensional collapsible
polyhedron (in particular, any cone over a graph) can be embedded
in a product of two trees. In section 2F we construct a
2-dimensional polyhedron embeddable in a product of two curves
which is not embeddable in any product of two graphs.

\ms \centerline{2A. {\it Definitions and general properties of
certain generalized manifolds}}\ms

Let $K$ be a $CW$ complex. Then the open cells of $K$ (that is,
the interiors $\kn\sigma$ of the cells $\sigma \in K$) form a {\it
partition} of $|K|$, i.e. they are mutually disjoint and cover
$|K|$. It follows from the definition of a $CW$ complex that for
each skeleton $K^{(n)}$ the space $|K^{(n)}|$ is the union of open
cells with dimension $\leq n$. If $\dim\sigma =n$ then
$\partial\sigma\subset |K^{(n-1)}|$ and
$\kn\sigma\cap|K^{(n-1)}|=\emptyset$. A cell $\sigma \in K$ is
said to be {\it proper} if it is a union of open cells of $K$. If
each cell of $K$ is proper then $K$ is said to {\it have proper
cells}. Two cells $\sigma, \tau\in K$ are said to be {\it
incident} if one is a subset of the other; if $\sigma\subset\tau$
then $\sigma$ is called a {\it face} of $\tau$.

\proclaim{Proposition 2A.1} A CW complex K has proper cells if and
only if for each two cells $\sigma, \tau \in K$ the condition
$\kn\sigma\cap\tau\neq\emptyset$ implies $\sigma\subset\tau$
$($that is, $\sigma$ is a face of $\tau$$)$. \qed\endproclaim

\proclaim{Corollary 2A.2} If $K_1,\cdots, K_n$ are CW complexes
with proper cells then $K_1\kw\cdots\kw{}K_n $ has proper cells as
well. \qed\endproclaim

Notice that there exist finite $CW$ complexes with proper cells
which are not regular. (The standard $CW$ structure on $\Sb^1$ has
this property.) However, the converse holds:

\proclaim{Lemma 2A.3} Any regular CW complex K has proper cells.
\endproclaim

 \demo{Proof} We must show that any cell $\tau\in K$
is a union of open cells. By induction we may assume that this
holds for all cells with dimension $<n+1=\dim \tau$, $n\geq 0$.
Note that $\partial\tau$ is a subset of the union of the $n$-cells
$\sigma_1,\cdots,\sigma_k$ of $K$ such that
$\kn\sigma_i\cap\tau\neq\emptyset$. It remains to show that
$\sigma_i\subset\tau$ for each $i$. To this end, fix $i$. As
$\partial\tau\cap\kn\sigma_i$ is a non-void closed subset of
$\kn\sigma_i$, to get the conclusion, it is enough to prove that
$\partial\tau\cap\kn\sigma_i$ is open in $\kn\sigma_i$. Since
$\partial\tau$ is an $n$-sphere in $|K^{(n)}|$ and $\kn\sigma_i$
is open in $|K^{(n)}|$, for each point
$x\in\partial\tau\cap\kn\sigma_i$ there is an open $n$-disc
containing $x$ and wholly lying in this intersection. By the
Brouwer Domain Invariance theorem such disc is a neighborhood of
$x$ in $\kn\sigma_i$. Thus $\partial\tau\cap\kn\sigma_i$ is open
in $\kn\sigma_i$, which completes the proof.\qed
\enddemo\ms

  From the classic Borsuk Separation Theorem relating closed sets separating
${\Sb}^n$ to their essential mappings into ${\Sb}^{n-1}$ (cf.
\cite{E-S, p. 302}) we infer the following fact.

\proclaim{Lemma 2A.4} For any $n$-manifold $M$ and a point $x_0
\in M$ there is an open neighborhood $V$ of $x_0$ in $M$ such that
every closed subset $F$ of $M$ separating $M$ between $x_0$ and $M
\setminus V$ admits an essential map into $\Sb^{n-1}$. $($In fact,
this holds for every neighborhood $V$ which is an open n-disc.$)$
\qed
\endproclaim

The revealed property can be used to define, for each natural
$n\ge 1$, a new class of $n$-dimensional continua (comprising in
particular all connected closed $n$-manifolds) playing an
important role in our investigations. Namely, an $n$-dimensional
continuum $X$ is said to be a {\it quasi} $n$-{\it manifold at a
point} $x \in X$ if there is an open neighborhood $V$ of $x$ in
$X$ such that every closed subset $F$ of $X$ with $\dim F \leq
{n-1}$ separating $X$ between $x$ and $X \setminus V$ admits an
essential map into $\Sb^{n-1}$. (Recall that a closed set
$F\subset X$ is said to {\it separate X between subsets} $A$ and
$B$ if there exist disjoint open sets $U$ and $V$ such that
$X\setminus F=U\cup V$, $A\subset U$ and $B\subset V$.) Notice
that $X$ is a quasi 1-manifold at $x\in X$ if and only if $x$ is
not an endpoint of $X$. (A point of a space is said to be its {\it
endpoint} if that point admits arbitrarily small open
neighborhoods whose boundaries are one-point sets.) If $X$ is a
quasi $n$-manifold at $x$ then every neighborhood of $x$ is
$n$-dimensional. If $V$ is as in the first definition then any
other open neighborhood $W$ of $x$ contained in $V$ has the
separation property as well. Also note that if a closed set
$F\subset V$ separates $\overline{V}$ between $x$ and $\partial V$
then $F$ separates $X$ between $x$ and $X\setminus V$. If $X$ is a
quasi $n$-manifold at every point of $X$ then it is called a {\it
quasi} $n$-{\it manifold}. Notice that an $n$-dimensional
continuum which is a union of quasi $n$-manifolds is a quasi
$n$-manifold as well.

An $n$-dimensional continuum $X$ is said to be a {\it para
n-manifold} if each point of $X$ belongs to an open $n$-disc lying
in $X$ (not necessarily open in $X$). In other words, $X$ is a
union of open $n$-discs.

If $X$ is an $n$-dimensional continuum and a point $x\in X$ is an
element of an open $n$-disc lying in $X$ then $X$ is a quasi
$n$-manifold at $x$. (This follows from the observation that if
$E$ is an open disc in $X$ then $E$ is open in $\overline{E}$.)
Hence any para $n$-manifold is a quasi $n$-manifold. This simple
criterion can be used to detect many interesting quasi
$n$-manifolds which are not $n$-manifolds. For instance, both the
$\lq \lq$Bing house" and the $\lq\lq$dunce hat" are para
2-manifolds, so they are also quasi 2-manifolds. These examples
are widely known primarily for being 2-dimensional contractible
and not collapsible polyhedra.

And it is convenient to introduce other generalizations of
$n$-manifolds. First, for an $n$-dimensional continuum $X$ define
the following subsets: \ms $M_n(X)=\{x \in X: \text {\it x is an
element of an open n-cell lying in X and open in X} \}$; \ms
$P_n(X)= \{x \in X: \text{\it x is an element of an open n-cell
lying in X} \}$.\ms

\noi Thus $M_n(X)=X$ if and only if $X$ is a closed $n$-manifold,
and $P_n(X)=X$ if and only if $X$ is a para $n$-manifold. The set
$M_n(X)$ is a maximal $n$-manifold with empty boundary lying in
$X$ as an open subset, we call it the $n$-{\it manifold part of
}$X$.

Now we define the generalizations. An $n$-dimensional continuum
$X$ is said to be a {\it pseudo $n$-manifold} ({\it ramified
$n$-manifold}, respectively) if $M_n(X)$ ($P_n(X)$, respectively)
is dense in $X$ and $\dim [X\setminus M_n(X)]\leq n-2$ ($\dim
[X\setminus P_n(X)]\leq n-2$, respectively). If, in addition,
$M_n(X)$ ($P_n(X)$, respectively) is connected then $X$ is said to
be a {\it simple pseudo n-manifold} ({\it simple ramified
n-manifold}, respectively). Clearly, every pseudo $n$-manifold is
a ramified $n$-manifold.

Let $K$ be a regular $CW$ complex, and let $n\ge1$. The set of
$n$-cells of $K$ is said to be {\it chain connected} if every two
$n$-cells of $K$ can be joined by a finite sequence of $n$-cells
of $K$ such that every two successive cells meet along an
$(n-1)$-cell (i.e. both are incident with such a cell). If each
$(n-1)$-cell of $K$ is a face of an $n$-cell, the above property
holds if and only if the space $|K^{(n)}|\setminus |K^{(n-2)}|$ is
connected.

(In the literature, the term {\it pseudo $n$-manifold} is often
used to mean a polyhedron which is a simple pseudo $n$-manifold,
see Proposition 2A.5(i) below.)

\proclaim{Proposition 2A.5} Let K be an n-dimensional finite
regular CW complex. Then

\ss {\rm(i)} $|K|$ is a pseudo $n$-manifold $($simple pseudo
$n$-manifold, {\it respectively}$)$ if and only if each cell of K
is a face of an n-cell, each $(n-1)$-cell of K is incident with
exactly two n-cells $($and the set of $n$-cells of $K$ is chain
connected, respectively$)$;\ss

{\rm(ii)} $|K|$ is a ramified $n$-manifold $($simple ramified
$n$-manifold, respectively$)$ if and only if each cell of K is a
face of an n-cell, each $(n-1)$-cell of K is incident with at
least two n-cells $($and the set of $n$-cells of $K$ is chain
connected, respectively$)$;\ss

{\rm(iii)} If $|K|$ is a quasi $n$-manifold then $|K|$ is a
ramified $n$-manifold.\qed\endproclaim

\proclaim {Corollary 2A.6} Let K be an n-dimensional finite
regular CW complex. If $|K|$ is a simple pseudo $n$-manifold and L
is a subcomplex of K such that $|L|$ is a ramified n-manifold then
$L=K$. In particular, the conclusion holds if $|K|$ is an
n-manifold.\qed\endproclaim

Our main observation in this section is Theorem~2A.10 below which
describes a basic property of quasi $n$-manifolds and ramified
$n$-manifolds lying in $n$-dimensional polyhedra. In its proof
Lemma 2A.8 is needed. In the proof of Lemma 2A.8 we need in turn
Lemma~2A.7 below which is relatively simple but not trivial. Lemma
2A.7 is certainly known to many topologists, and can be proved
using various arguments. We supply possibly the shortest one.

\proclaim {Lemma 2A.7} No proper closed subset of $\Sb^n$ admits
an essential map into $\Sb^n$.\endproclaim

\demo{Proof} We may assume that $n\ge1$. Consider a proper closed
subset $F$ of $\Sb^n$. Since $\Sb^n$ is a compact subset of
$\R^{n+1}$, the set $F$ can be regarded as a compact subset of
$\R^{n+1}$. Then note that $\R^{n+1}\setminus F$ is connected. So,
by the Borsuk Separation Theorem \cite{E-S, p. 302}, $F$ admits no
essential map into $\Sb^n$.\qed
\enddemo

\proclaim{Lemma 2A.8}{\rm(a)} Let $X$ be a quasi $n$-manifold at a
point $x$. If $U$ is a neighborhood of $x$ in $X$ and $h:U\to
\R^n$ is an embedding, then $h(U)$ is a neighborhood of $h(x)$ in
$\R^n$.

{\rm(b)} Let $X$ be a ramified n-manifold and let U be a non-void
open subset of X. If $h:U\to \R^n$ is an embedding such that
$h(U)$ is closed, then $h(U)=\R^n$.
\endproclaim

\demo{Proof of $(${\rm a}$)$} Suppose, to the contrary, that
$h(x)\in\partial h(U)$. Let $V$ be an open neighborhood of $x$ in
$X$ with $\overline{V}\subset U$ such that any closed subset $F$
of $X$ separating $X$ between $x$ and $X\setminus V$ admits an
essential map into $\Sb^{n-1}$. Since $h(U\setminus V)$ is a
closed subset of $h(U)$ not containing $h(x)$ there is an open
ball $B(h(x),\varepsilon_0)$ in $\R^n$ such that
$B(h(x),\varepsilon_0)\cap h(U\setminus V)=\emptyset$. Also, there
is a sphere $S=\partial B(h(x),\varepsilon)$, $0 < \varepsilon
<\varepsilon_0$, such that $S\nsubseteq h(U)$. Since $S$ separates
$\R^n$ between $x$ and $h(U\setminus V)$ , the intersection
$E=S\cap h(U)$ separates $h(U)$ between $h(x)$ and $h(U\setminus
V)$. Note that $E$ is a proper subset of $S$. As $E=S\cap
h(V)=S\cap h(\overline V)$, the set $E$ is a compact subset of
$h(V)$ which separates $h(\overline V)$ between $h(x)$ and
$h(\partial V)$. It follows that $F=h^{-1}(E)$ is a compact subset
of $V$ which separates $\overline V$ between $x$ and $\partial V$.
Hence $F(\subset V)$ is a closed subset of $X$ which separates $X$
between $x$ and $X\setminus V$. By our choice of $V$, there is an
essential map $F\to \Sb^{n-1}$. But $F(\approx E)$ is homeomorphic
to a proper closed subset of $\Sb^{n-1}$, which contradicts Lemma
2A.7.\enddemo

\demo{Proof of $(${\rm b}$)$} It is enough to show that $h(U)$ is
dense in $\R^n$. Let $V$ = int $h(U)$. First we show that $V$ is
dense in $h(U)$. In fact, since $R(X)$ is dense $X$ and $U$ is
open, the image $h(R(X)\cap U)$ is dense in $h(U)$. On the other
hand, $h(R(X)\cap U)$ is open in $\R^n$, by the Brouwer Invariance
of Domain Theorem. Hence $h(R(X)\cap U)\subset V$, so $V$ is dense
in $h(U)$. Therefore, $\partial V=h(U)\setminus V$. Consequently,
$\partial V \subset h(U\setminus R(X))$, hence $\dim
\partial V \leq n-2$, by the definition of a ramified
$n$-manifold. Therefore, $\R^n \setminus \partial V$ is connected
(cf. \cite {E, Theorem 1.8.13, p. 77}). It follows that $V$ is
dense in $\R^n$, for otherwise $\partial V$ separates $\R^n$.
Hence $h(U)=\R^n$, which completes the proof. \qed\enddemo

A modification of the above argument gives the following

\proclaim{Corollary 2A.9} If $f:X\to Y$ is an embedding of a
ramified n-manifold into a simple pseudo n-manifold, then
$h(X)=Y$.\qed\endproclaim

\proclaim{Theorem 2A.10} Let $X$ be either a quasi $n$-manifold or
a ramified $n$-manifold. If $f:X\to |K|$ is an embedding, where
$|K|$ is an $n$-dimensional polyhedron, then $f(X)=|L|$, where $L$
is a subcomplex of $K$.
\endproclaim

\demo{Proof} It suffices to prove that $f(X)$ is a union of
$n$-cells of $K$. To this end, it is enough to show that each
point $y\in f(X)$ is an element of an $n$-cell of $K$ which lies
in $f(X)$. Notice that there is an open neighborhood $V$ of $y$ in
$f(X)$ such that, for each cell $\sigma \in K$, the condition
$V\cap \sigma\neq\emptyset$ implies $y \in \sigma$. It follows
from our hypothesis that ${\dim}V=n$. Hence there is an $n$-cell
$\sigma_0\in K$ such that $V\cap\kn\sigma_0\neq\emptyset$. Hence
$y\in \sigma_0$. It remains to show that $\sigma_0\subset f(X)$.
As $f(X)$ is closed, it is enough to show that $\kn\sigma_0\subset
f(X)$. Notice that $f(X) \cap \kn \sigma_0$ is a non-void closed
subset of $\kn\sigma_0$. On the other hand, this set is open in
$f(X)$. It follows that $U=f^{-1}(\kn\sigma_0)$ is non-void and
open in $X$, and $f(U)=f(X)\cap\kn\sigma_0$. Since $\kn\sigma_0
\approx \R^n$, by Lemma 2A.8 we infer that $f(X)\cap\kn\sigma_0$
is also open in $\kn\sigma_0$, in both cases under discussion.
Consequently, $f(X) \cap \kn\sigma_0=\kn\sigma_0$. Hence
$\kn\sigma_0 \subset f(X)$, as desired.\qed\enddemo

We have already noted that any polyhedron which is a quasi
$n$-manifold is a ramified $n$-manifold as well. The converse
holds for $n=1,2$ but in general fails: the suspension of the
"dunce hat" (or the "Bing house")is a polyhedron which is a
ramified $3$-manifold but is not quasi $3$-manifold.

Let $|K|$ be a ramified $n$-manifold. By its {\it combinatorial
component} we mean any maximal simple ramified $n$-manifold in
$|K|$. One easily sees that $|K|$ is the union of its
combinatorial components, and any two different combinatorial
components meet in a subpolyhedron of dimension $\leq{n-2}$.

A connected graph $P$ is a ramified 1-manifold if and only if it
has no endpoint. In such a case $\rank H^1(P)\ge1$, and $\rank
H^1(P)=1$ if and only if $P$ is a circle. For any $n\ge1$ there
exist only finitely many topological types of ramified 1-manifolds
$P$ with $\rank H^1(P)=n$. For any connected graph $P$ we have
$H^1(P)=n$ if and only if $P$ is homotopy equivalent to a bouquet
of $n$ circles. Any graph which is a ramified 1-manifold is
simple.

We have the following diagram of subclasses of all connected
polyhedra:
$$
\CD \{para \ n{-}manifolds\}\ & \supset & \ \ \{simple\
para \ n{-}manifolds\} \\
    \cap & {} & {}\\
\{quasi \ n{-}manifolds\} \ & {} &  \cap \\
    \cap & {} & {}\\
\{ramified \ n{-}manifolds\}\ & \supset & \ \ \{simple\
ramified \ n{-}manifolds\} \\
\cup & {} & \cup \\
\{pseudo \ n{-}manifolds\} & \supset & \ \ \{simple\ pseudo\
n{-}manifolds\}.
\endCD
$$

\ms\centerline{2B. {\it On locally connected quasi manifolds.}}
\centerline{\it From embeddings into products of curves to
embeddings into products of graphs}\ms

 Here we prove a useful theorem on factorization of embeddings of quasi
manifolds into products of curves through embeddings into product
of graphs.

\proclaim{Theorem 2B.1} Let $X$ be a locally connected quasi
$n$-manifold such that $H^1(X)$ has finite rank. If $f = (f_1,
\cdots ,f_n):X \to Y_1 \times \cdots \times Y_n$ is an embedding
of $X$ in the product of $n$ curves, then there exist mappings $g
= (g_1, \cdots ,g_n):X \to P_1 \times \cdots \times P_n$ and $h=
h_1 \times \cdots \times h_n:P_1 \times \cdots \times P_n \to Y_1
\times \cdots \times Y_n$ such that $f_i = h_i \circ g_i$ for each
$i=1,\cdots ,n$ $($hence $f=h \circ g$$)$, where $g_i:X \to P_i$
is a monotone surjection, $P_i$ is a graph with no endpoint
$($that is, $P_i$ is a quasi $1$-manifold$)$, and $h_i:P_i\to Y_i$
is $0$-dimensional.

In particular, if $X$ is embeddable in a product of $n$ curves,
then there exists an embedding $(g_1, \cdots ,g_n):X \to P_1
\times \cdots \times P_n$, where each $g_i:X \to P_i$ is a
monotone surjection, $P_i$ is a graph with no endpoint, and $
\rank$ $H^1(P_i)\le \rank$ $H^1(X)$.
\endproclaim

\ms\noindent{\bf Note.} It follows that {\it if} $f_i:X \to Y_i$
{\it is monotone, then} $f_i(X)$ {\it is a graph.} In fact, in
this case $f_i(X)=h_i(P_i)$ and $h_i:P_i\to Y_i$ is an embedding.
If $f_i$ is not monotone then $f_i(X)$ need not be a graph, see
Example 2B.3. The proof given below shows that {\it if $X$ is a
non-degenerate connected polyhedron $($or any non-degenerate
locally connected continuum whose $H^1(X)$ has finite rank$)$ then
$f_i(X)$ is a local dendrite.} \qed\bs

\demo{Proof} By the Whyburn factorization theorem, there is a
factorization $f_i=h_i \circ g_i$ ,
$$
X @>{g_i}>> P_i @>{h_i}>> Y_i,
$$
where $g_i$ is a monotone surjection, and $h_i$ is
$0$-dimensional. Since $Y_i$ is 1-dimensional and $h_i$ is
0-dimensional we infer that $\dim P_i \leq1$ (by a theorem of
Hurewicz). Clearly, $g = (g_1, \cdots ,g_n):X \to P_1 \times
\cdots \times P_n$ is an embedding. Since $\dim X=n$, it follows
that $\dim P_i>0$ for each $i$. Therefore, $P_i$ is a locally
connected curve, as $g_i$ is a surjection. Since $g_i$ is a
monotone surjection and $H^1(X)$ has finite rank, $P_i$ is
actually a local dendrite (see \cite{Kr, Lemma 3.1}). Hence each
point of $P_i$ has a closed neighborhood which is a dendrite.
First we show that\ms

 (*) {\it $P_i$ has no endpoint}.\ms

\noi For suppose $P_i$ has an endpoint $z_0$. Since $g_i(X) =
P_i$, there is a point $x_0 \in X$ such that $g_i(x_0) = z_0$.
Since $X$ is a quasi $n$-manifold at $x_0$ there is an open
neighborhood $V$ of $x_0$ in $X$ such that

$(1)$ {\it any closed $(n-1)$-dimensional subset of} $X$ {\it
separating} $X$ {\it between} $x_0$ {\it and} $X \setminus V$ {\it
admits an essential map to} $\Sb^{n-1}$.

Now we shall show that there is an open neighborhood $U$ of
$g(x_0)$ in $P_1 \times \cdots\times P_n$ such that \roster \item
"(2)" $\overline{U} \cap g(X \setminus V)=\emptyset$, \item "(3)"
$\dim \partial U = n-1$, \item "(4)" $\partial U$ {\it is
contractible}.
\endroster

\noi To construct $U$ we assume, without loss of generality, that
$i=1$. Then $g(x_0)=(y_1,y_2,\cdots,y_n)$, where $y_1=z_0$. Note
that $g(V)$ is a neighborhood of $g(x_0)$ in $g(X)$, hence any
small enough $U$ satisfies (2). Since $z_0$ is an endpoint of
$P_1$, and each $P_j$ is a local dendrite, there exist sets
$U_1,\cdots,U_n$, as small as we please, such that each $U_j$ is
an open and connected neighborhood of $y_j$ in $P_j$ with
$\partial U_j$ finite, each $\overline{U_j}$ is a dendrite, and
$\partial U_1$ is a one-point set. Then the set
$U=U_1\times\cdots\times U_n$ has the desired properties. In fact,
as $U_j$'s are small, $U$ satisfies (2). Then note that
$$\partial U=\bigcup_{j=1}^{n} (
\overline{U_1}\times\cdots \times \overline{U_{j-1}} \times
\partial U_j\times \overline{U_{j+1}}\times\cdots\times
\overline{U_n}).$$ \noi Hence (3) follows. To prove (4), note that
$(\partial U_1)\times\overline{U_2}\times\cdots\times
\overline{U_n}$ is a strong deformation retract of $\partial U$
(because $\partial U_1$, as a one-point set, is a strong
deformation retract of $\overline{U_1}$). Hence (4) follows from
the fact that $(\partial
U_1)\times\overline{U_2}\times\cdots\times \overline{U_n}$ is
contractible.

Now consider the set $F =\partial_{g(X)} (U \cap g(X))$. Observe
that it is a closed subset of $g(X)$ such that \roster \item "(5)"
$F \subset \partial U$, \item "(6)" $F$ {\it separates} $g(X)$
{\it between} $g(x_0)$ {\it and} $g(X \setminus V)$.
\endroster
It follows that \roster \item "(7)" $g^{-1}(F)$ {\it is closed},
$(n-1)$-{\it dimensional, and separates X between} $x_0$ {\it and}
$X\setminus V$.\endroster

Now we are ready to complete the proof of (*). Note that by (1)
and (7) there is an essential map $\varphi: F \to \Sb^{n-1}$. By
(3) and (5) there is a continuous extension $\varphi^*:
\partial U \to \Sb^{n-1}$ of $\varphi$. However, by (4),
$\varphi^*$ is null-homotopic, hence so is $\varphi$, a
contradiction. This proves (*).

Next we show that\ms

(**) {\it $P_i$ is a graph}.

\ms To prove (**) recall that $P_i$ is a local dendrite. Since
$P_i$ has no endpoint, it contains a circle. (Otherwise it is a
dendrite, hence contains an endpoint.) It follows that the union
of all simple closed curves in $P_i$ is a (not necessarily
connected) graph. Enlarging this set by the union of a finite
collection of arcs (e.g., adding arcs in $P_i$ irreducibly
connecting different components of the union), we get a connected
graph $Q_i (\subset P_i)$ such that for each component $C$ of $P_i
\setminus Q_i$ we have

$(8)$ $\overline{C}$ {\it is a dendrite and} $\partial C$ {\it
consists of a single point}.

\noi To prove (**) it suffices to show that $Q_i = P_i$.

Suppose, on the contrary, that $P_i \setminus Q_i \neq \emptyset$.
Then consider a component $C$ of $P_i \setminus Q_i$. It is an
open set in $P_i$. By (1), $\overline{C}$ is a dendrite and
$\partial C$ is a one-point set. There is a point $z_0 \in C$
which is an endpoint of the dendrite $\overline{C}$. Since $C$ is
open, $z_0$ is an endpoint of $P_i$ as well. This contradicts (*)
and ends the proof of (**).

As $g_i$ is a monotone surjection, the induced homomorphism
$H^1(g_i):H^1(P_i)\to H^1(X)$ is a monomorphism by the
Vietoris-Begle Theorem for $n=1$ (see e.g. \cite{Sp, 6.9, Theorem
15}). Therefore, rank $H^1(P_i) \le$ rank $H^1(X)$. This completes
the proof. \qed
\enddemo

It is well known that for any closed $n$-manifold $M$ the group
$H^1(M)$ has finite rank. Consequently, Theorem 2B.1 implies the
following

\proclaim{Corollary 2B.2} If a closed $n$-manifold is embeddable
in a product of $n$ curves, then it is also embeddable in a
product of $n$ graphs.\qed
\endproclaim

\noi {\bf{Note.}} This corollary answers a question posed (for
surfaces) by R. Cauty \cite{C}.\qed\ms

\proclaim{Example 2B.3} There exist a curve $X$ which is not a
graph, and a closed orientable surface $M$ in the product $X
\times X$ such that both projections $pr_i:X \times X \to X$ map
$M$ onto $X$. Moreover, $M$ is invariant under the canonical
involution on $X \times X$ which interchanges the coordinates.
\endproclaim

\demo{Proof} First we construct a closed orientable surface $N$ in
the product $Y_1 \times Y_2$ of two curves such that \roster \item
"(1)" $Y_1$ {\it is not a graph}, \item "(2)" {\it the
projections} $q_i:Y_1\times Y_2 \to Y_i$ {\it map} $N$ {\it onto}
$Y_i$.
\endroster
Define $Y_1$ to be the union $Y_1=\alpha_0 \cup \alpha_1
\cup\beta_1 \cup \beta_2$ of four arcs with common endpoints $a,b$
such that $\alpha_0 \cup \alpha_1 \cup \beta_1$ ia a
$\theta$-curve, $\beta_2\cap(\kn\alpha_0 \cup
\kn\alpha_1)=\emptyset$, and $\beta_1\cap \kn\beta_2$ is a compact
set with infinitely many components. Then $Y_1$ satisfies (1).
Define $Y_2$ to be a graph given by the formula: $Y_2=T_0 \cup
T_1$, where $T_0$, $T_1$ are two oriented circles whose
intersection $T_0 \cap T_1=L_0 \cup L_1$, where $L_0,L_1$ are
disjoint oriented arcs coherently oriented with each $T_i$. In
such a case, each $T_i$ can be presented as the union of four arcs
with disjoint interiors, $T_i=A_i \cup B_i \cup L_0 \cup L_1$,
such that $S_1=A_0 \cup A_1$ and $S_2=B_0 \cup B_1$ are disjoint
circles. We define the surface $N$ in $Y_1 \times Y_2$ by the
formula
$$
N=\alpha_0 \times T_0 \cup \alpha_1 \times T_1 \cup \beta_1\times
S_1 \cup \beta_2\times S_2.
$$

\noindent One easily verifies that $N$ is an orientable surface
satisfying $(2)$.

To construct the promised example we proceed as follows. Choose a
homeomorphism $h:\alpha_0 \to A_0$. Then define $X$ to be the
quotient space $X=(Y_1\sqcup Y_2)/x\sim h(x)$ for each $x\in
\alpha _0$). By $(1)$ we infer that $X$ is a curve but not a
graph. Let $X_i = h_i(Y_i)$, where $h_i:Y_i \to X$ are canonical
embeddings. Clearly, $X = X_1 \cup X_2$ and $X_1 \cap X_2 = A$,
where $A = h_1(\alpha_0) = h_2(A_0)$ is an arc. Let $t:Y_1 \times
Y_2 \to Y_2 \times Y_1$ denote the map given by $t(y,z) = (z,y)$.
Then we define $M (\subset X \times X)$ as follows:
$$
M = [(h_1 \times h_2)(N) \cup (h_2 \times h_1)(t(N))] \setminus
(\kn A \times \kn A).
$$
One easily verifies that $M$ is invariant under canonical
involution on $X \times X$. Since $\alpha_0 \times A_0 \subset N
\subset Y_1 \times Y_2$, $A_0 \times \alpha_0 \subset t(N) \subset
Y_2 \times Y_1$. Hence $A \times A = (h_1 \times h_2)(\alpha_0
\times A_0) \subset (h_1 \times h_2)(N) \subset (h_1 \times
h_2)(Y_1 \times Y_2) = X_1 \times X_2$. Likewise, $A \times A
\subset (h_2 \times h_1)(t(N)) \subset X_2 \times X_1$. Moreover,
$(X_1 \times X_2) \cap (X_2 \times  X_1) = A \times A$. Hence we
have
$$
(h_1 \times h_2)(N) \cap (h_2 \times h_1)(t(N)) = A \times A.
$$
Thus $M$ is the connected sum of orientable surfaces $(h_1 \times
h_2)(N)$ and $(h_2 \times h_1)(t(N))$, hence it is an orientable
surface. Applying $(2)$, we easily see that both projections
$pr_i:X \times X \to X$ map $M$ onto $X$. \qed
\enddemo

\ms In connection with the above proof let us notice the following
fact.

\ms \noi{\bf Note.} Let $M$ be any compactum lying in the product
$P_1\times P_2$ of two graphs. If $A$ is an arc in $P_1$ with
$p_1^{-1}(A)=A\times(S_1\cup\dots\cup S_k)$, $k\ge2$, where
$S_1,\dots,S_k$ are disjoint circles in $P_2$, then $M$ can be
embedded in the product $P_1'\times P_2$ in such a way that $P_1'$
is a curve and the image of $M$ under the projection $P_1'\times
P_2\to P_1'$ is not a graph. (In fact, $P_1'$ can be obtained from
$P_1$ by adding an arc $B$ with the same endpoints as that of $A$
in such a way that $A\cup B$ is not a graph.) \qed\ms

The surface $M$ constructed in Example 2B.3 meets the diagonal of
$X\times X$. Below we present another example of a surface in the
product $P\times P$, where $P$ is a graph, which is disjoint from
the diagonal of the product and is invariant under the canonical
involution on $P\times P$.

\proclaim{Example 2B.4} There exist a graph $P$ and a closed
orientable surface $M$ in $P \times P$ such that: $M$ is disjoint
with the diagonal of $P\times P$, both projections $pr_i:P \times
P \to P$ map $M$ onto $P$, and $M$ is invariant under the
canonical involution on $P \times P$.
\endproclaim

\demo{Proof} Fix any number $n\ge 4$. The graph $P$ is defined to
be a subset of $\Sb^1\times I$ given by
$$P=(\Sb^1\times \{0,1\}) \cup \{z_0,\cdots,z_{n-1}\}\times I,$$
where $z_j=\exp(2\pi i \frac{j}{n})$, $j=0, \cdots, n-1$. Then
define arcs $A_j\times\{0\}$, $A_j\times\{1\}$, $I_j$ and circles
$S_j$ in $P$ as follows:
$$A_j=\{\exp(2\pi it):t\in [\frac{j}{n},\frac{j+1}{n}]\}, \;\;
I_j=\{z_j\}\times I,\;\; S_j=I_j \cup (A_j \times \{0,1\}) \cup
I_{j+1}.$$ \noi (All indices in this construction are reduced
modulo $n$.) Finally, define tori $T_j$ to be the subsets of
$P\times P$ given by
$$T_j=S_j\times S_{j+2}.$$ \noi Notice that the intersection
$$D_j=T_i \cap T_{j+1}=I_{j+1}\times I_{j+3}$$
is a disc. Now we are ready to define the surface $M$, put
$$M=(T_0 \cup \cdots \cup T_{n-1})\setminus (\kn D_0\cup\cdots\cup\kn D_{n-1}).$$
One easily verifies that $M$ has all the desired properties.\qed
\enddemo\ms

\centerline{2C. {\it Ramified manifolds in products of graphs}}\ms

Here we establish some properties of ramified $n$-manifolds lying
in products of $n$ graphs. These properties will find essential
applications in the next section and in a subsequent paper.\ms

\centerline{\it{Throughout this section we consider fixed graphs
$P_1=|K_1|$, $\dots$,$P_n =|K_n|$, $n\ge2$,}} \centerline
{\it{where each $K_i$ is a regular $1$-dimensional $CW$
complex.}}\ms \noi [So, each 1-cell of $K_i$ (an {\it edge}), is
an arc; its endpoints are called {\it vertices}.] By
$K_1\kw\cdots\kw{}K_n$ we denote the $CW$ cell structure on
$P_1\times\cdots\times{}P_n$ defined by
$$K_1\kw\cdots\kw{}K_n=\{\sigma_1\times\cdots\times\sigma_n:
\sigma\in{}K_1,\cdots,\sigma_n\in{}K_n \}. $$
\centerline{\it{Also, we consider a fixed ramified $n$-manifold
$M=|K(M)|$ lying in}}
\centerline{\it{$P_1\times\cdots\times{}P_n$, where $K(M)$ is a
subcomplex of $K_1\kw\cdots\kw{}K_n$.}}\ms \noi By Proposition 2A.
5(ii) we have

\proclaim{Property (a)} $M$ is the union of $n$-cells of $K(M)$.
\qed\endproclaim

We adopt the following notation. For a {\it non-void subset} $J$
of the index set $\{1,\cdots,n\}$ let: \ss

- $P_J=\prod _{j\in J}P_j$,

- $K_J$ denote the cell structure on $P_J$ induced by
$\{K_j:j\in{}J\}$,

- $p_J$ denote the restriction to $M$ of the projection
$pr_J:P_1\times\cdots\times{}P_n\to{}P_J$ (in particular,
$p_{\{1,\cdots,n\}}:M\to{}P_1\times\cdots\times{}P_n$ is the
inclusion mapping,

- $n_J=|J|$,

- $J^c=\{1,\cdots,n\}\setminus{}J$ (therefore, $n_{J^c}=n-n_J$).

\noindent Notice that $P_{\{j\}}=P_j$, $K_{\{j\}}=K_j$ and
$pr_{\{j\}}=pr_j$. We abbreviate $p_{\{j\}}$ to $p_j$.

\ms For any cell
$\sigma=\sigma_1\times\cdots\times\sigma_n\in{}K$,
$\sigma_j\in{}K_j$, the restriction $p_J|\sigma$ is the projection
onto $\sigma_J=\prod_{j\in{}J}\sigma_j$. In this sense we say that
$p_J$ {\it preserves} the cell structures $K$ and $K_J$. It
follows from Property (a) that

\proclaim{Property (b)} $p_J(M)=|K'_J|$, where $K'_J$ is a
subcomplex of $K_J$. Moreover, $|K'_J|$ is a ramified
$n_J$-manifold; if $M$ is a $($simple$)$ pseudo $n$-manifold, then
$p_J(M)$ is a $($simple$)$ pseudo $n_J$-manifold. {\rm (In the
sequel we abbreviate $K'_{\{j\}}$ to $K'_j$.)}\qed\endproclaim

\ms From this point on to the end of this section we assume that
$J$ {\it is a proper non-void subset of} $\{1,\cdots,n\}$.

\ms For every cell $\tau\in{}K'_{J^c}$, we define $P_{J}(\tau)$ to
be the union of all $n_{J}$-cells $\sigma\in{K}_{J}$ such that
$\sigma\times\tau\subset{}M$. From Property (a) we infer that

\proclaim{Property (c)}
$M=\bigcup\{P_{J}(\tau)\times\tau:\tau\in{}K'_{J^c}\ is\ an\
n_{J^c}{-}cell\}$. \qed\endproclaim

\proclaim{Property (d)} If $\tau$ is a face of a cell
$\tau'\in{}K'_{J^c}$ then $P_{J}(\tau)\supset{}P_{J}(\tau')$.
\qed\endproclaim

\ms In addition, for any point $y\in{p}_{J^c}(M)$, we define
$P_{J}(y)$ to be the set $P_{J}(y)=\{x\in{}P_{J}:(x,y)\in{}M\}$.
Thus, $P_{J}(y)\times\{y\}=p_{J^c}^{-1}(y)$, and
$P_{J}(y)\subset{}p_{J}(M)$.

\proclaim{Property (e)} For any cell $\tau\in K'_{J^c}$ and any
point $y\in\kn\tau$ we have $P_{J}(y)=P_{J}(\tau)=\bigcup
\{P_{J}(\tau'):\tau'\in{}K'_{J^c}\ is\ an\ n_{J^c}{-}cell\ with\
face\ \tau\}.$
\endproclaim

\demo{Proof} Note that
$P_{J}(y)\supset{}P_{J}(\tau)\supset{}P_{J}(\tau')$ for each cell
$\tau'\in{}K'_{J^c}$ with face $\tau$. So, it remains to justify
the inclusion
$$P_{J}(y)\subset\bigcup\{P_{J}(\tau'):\tau'\in{}K'_{J^c}\ is\
an\ n_{J^c}{-}cell\ with\ face\ \tau \}.$$ To this end, consider a
point $x\in{}P_{J}(y)$. Then $(x,y)\in{}M$. By Property (a),
$(x,y)$ belongs to an $n$-cell $\sigma\times\tau'\subset{}M$,
where $\sigma$ is an $n_{J}$-cell in $K_{J}$ and $\tau'$ is an
$n_{J^c}$-cell in $K'_{J^c}$. As $y\in\kn\tau$, $\tau$ is a face
of $\tau'$ by Proposition 2A.1, because $K_J$ has proper cells
(see Corollary 2A.2). It follows that
$x\in\bigcup\{P_{J}(\tau'):\tau'\in{}K'_{J^c}\ is\ an\
n_{J^c}{-}cell\ with\ face\ \tau \}$, which ends the proof. \qed
\enddemo

\proclaim{Property (f)}  For any cell $\tau\in{}K'_{J^c}$ the set
$P_{J}(\tau)$ is a finite disjoint union of ramified
$n_{J}$-manifolds in $P_{J}=|K_{J}|$. Moreover, if $M$ is a pseudo
$n$-manifold and $\tau$ is an $n_{J^c}$-cell then $P_{J}(\tau)$ is
a finite union of disjoint pseudo $n_{J}$-manifolds.
\endproclaim

\demo{Proof} If $\tau$ is an $n_{J^c}$-cell then both assertions
follow from the fact that each $(n{-}1)$-cell
$\sigma\times\tau\in{}K(M)$, where $\sigma\in{}K_{J}$ is an
$(n_{J}{-}1)$-cell, is a face of at least two (exactly two if $M$
is a pseudo $n$-manifold) $n$-cells $\sigma_1\times\tau,\
\sigma_2\times \tau\in{}K(M)$. Consequently, the first assertion
for arbitrary $\tau\in{}K'_{J^c}$ follows from Property (e). \qed
\enddemo

\proclaim{Property (g)} $p_{J}(M)=
\bigcup\{P_{J}(\tau):\tau\in{}K'_{J^c}\ is\ a\ k{-}cell\}$ for
each $k=0,\cdots,n_{J^c}$.
\endproclaim

\demo{Proof} For $k=n_{J^c}$ this follows from Property (c).
Applying Property (e) we obtain the general case. \qed\enddemo

\proclaim{Property (h)} If $P_{J}(w)$ is an $n_{J}$-manifold for
each vertex $w\in{}K'_{J^c}$, then $p_{J}(M)=P_{J}(w_0)$ for any
vertex $w_0$ of $K'_{J^c}$.
\endproclaim

\demo{Proof} By Property (g) (with $k=0$), it suffices to prove
that

\ss ($\ast$) $P_{J}(w_1)=P_{J}(w_2)$ for any two vertices
$w_1,w_2\in{}K'_{J^c}$. \ss

\noindent To this end, consider a 1-cell $\tau\in{}K'_{J^c}$ with
vertices $w$ and $w'$. Then, by Properties (f) and (d),
$P_{J}(\tau)$ is a finite union of ramified $n_{J}$-manifolds
contained in both $n_{J}$-manifolds $P_{J}(w)$ and $P_{J}(w')$. It
follows that $P_{J}(w)=P_{J}(\tau)=P_{J}(w')$, as no proper
ramified $n_{J}$-manifold  is contained in an $n_{J}$-manifold
(cf. Corollary 2A.6). Thus the condition ($\ast$) is a consequence
of the connectivity of the complex $K'_{J^c}$. \qed
\enddemo

\proclaim{Property (i)} If $p_{J}(M)$ is a simple pseudo
$n_{J}$-manifold then $M=p_{J}(M)\times{}p_{J^c}(M)$.
\endproclaim

\demo{Proof} For every $\tau\in{}K'_{J^c}$ the set $P_{J}(\tau)$
is a ramified $n_{J}$-manifold contained in the simple pseudo
$n_{J}$-manifold $p_{J}(M)$, so $P_{J}(\tau)=p_{J}(M)$. Hence the
assertion follows from Property(c).  \qed
\enddemo

\proclaim{Property (j)} Let $j\in \{1,\cdots, n\}$. Then $p_j(M)$
is a circle if and only if $P_j(v)$ is a circle for each vertex
$v\in K_{\{j\}^c}'$.
\endproclaim

\demo{Proof}This follows from Property (h) combined with Property
(i).\qed\enddemo\ms

Put
$$J_M=\{j\in\{1,\cdots,n\}:p_j(M)~ \text{\it {is a circle}}\}.$$
If $J_M$ is non-void then $M$ is said to {\it have projections
onto circles}, and any $p_j$, for $j\in J_M$, is said to be {\it a
projection of M onto a circle}.

\proclaim{Property (k)} If $J_M=\{1,\cdots,n\}$ then
$M=p_1(M)\times \cdots \times p_n(M)$ is an $n$-torus. If $J_M$ is
a proper non-void subset of $\{1,\cdots,n\}$ then
$M=p_{J_M}(M)\times p_{J_M^c}(M)$, where $p_{J_M}(M)=\prod_{j\in
J_M} p_j(M))$ is an $n_{J_M}$-torus, and $p_{J_M^c}(M)$
$($$\subset \prod_{j\in J_M^c}P_j$$)$ is a ramified
$n_{J_M^c}$-manifold which has no projection onto a circle.
\endproclaim

\demo{Proof} First notice that $p_{J_M}(M)$ is an $n_{J_M}$-torus
$\prod_{j\in J_M}p_j(M)$. This follows from Corollary 2A.6 because
$p_{J_M}(M)$ is a ramified $n_{J_M}$-manifold (see Property (b))
lying in the torus $\prod_{j\in J_M}p_j(M)$ of the same dimension.
Next notice that $M=p_{J_M}(M)\times p_{J_M^c}(M)$ by Property
(i). This ends the argument.\qed
\enddemo

\ms\centerline{2D. {\it Product structure of generalized manifolds
lying in products of curves}}\ms

The proof of the main result of this section, Theorem 2D.6, will
be preceded by a series of auxiliary lemmas.

\proclaim{Lemma 2D.1} Let M be a non-void compact subset of the
product $P\times Q$ and let $u:P\to Y$, $v:Q\to Z$ be mappings
such that $u\times v:P\times Q\to Y\times Z$ is injective on M.
Then we have:

\ss{\rm(i)} If every pair of sets $pr_Q(M\cap pr_P^{-1}(x))$, with
$x\in pr_P(M)$, has non-void intersection then $u$ is injective on
$pr_P(M)$ $;$

\ss{\rm(ii)} If $P$ is a graph, $Q$ is a finite product of graphs,
and $M$ is a ramified $n$-manifold then $pr_P(M)$ is a circle if
and only if $pr_Y((u\times v)(M))$ is a circle.
\endproclaim

\demo{Proof of {\rm(i)}} Suppose $u(x)=u(x')$ for some $x,x'\in
pr_P(M)$. It suffices to show that $x=x'$. By the hypothesis there
is $y\in pr_Q(M\cap pr_P^{-1}(x))\cap pr_Q(M\cap pr_P^{-1}(x'))$.
Hence $(x,y), (x',y)\in M$ and $(u\times v)(x, y)=(u\times v)(x',
y)$. Since $u\times v$ is injective on $M$ the conclusion follows.

\noi{\it Proof of {\rm(ii)}.} We consider all $P$, $pr_P(M)$, $Q$,
$pr_Q(M)$, and $M$ with natural $CW$ product structures as in
section 2C. First assume $pr_P(M)$ is a circle. We must show that
$pr_Y((u\times v)(M))$ is a circle as well. To this end, observe
that for any vertex $w\in pr_Q(M)$ the set $P(w)$ is a non-void
union of ramified 1-manifolds, see Property (f). Since $P(w)\times
\{w\}\subset M$, we also have $P(w)\subset pr_P(M)$. It follows
from our assumption that $P(w)=pr_P(M)$. Hence $P(w)$ is a circle
and $pr_P(M)=pr_P(P(w)\times \{w\})$. So, we have $pr_Y((u\times
v)(M))=u(pr_P(M))=u(pr_P(P(w)\times \{w\}))=pr_Y((u\times
v)(P(w)\times \{w\}))=u(P(w))$. Since $u\times v$ is injective on
$P(w)\times \{w\}$ and $pr_Y$ is injective on $u(P(w))\times
v(\{w\})$ it follows that $u(P(w))$ is a circle. This proves the
first implication.

Next assume $pr_Y((u\times v)(M))$ is a circle. We have to show
that $pr_P(M)$ is also a circle. By Property (h) it suffices to
show that $P(w)$ is a circle for each vertex $w\in pr_Q(M)$. But
this can be achieved using an argument analogous to that from the
above proof. This ends the proof of the second
implication.\qed\enddemo

\proclaim{Lemma 2D.2}Let $X$ be a compactum and let $A$ be a
closed subset of $X$. If $\dim X\leq m$ and $H^m(X)=0$ then
$H^n(A)=0$ for each $n\geq m$.
\endproclaim

\demo{Proof} Since $H^n(X)=0$ and $H^{n+1}(X, A)=0$, the
conclusion follows from the cohomology exact sequence of the pair
$(X, A)$. \qed\enddemo

\proclaim{Lemma 2D.3} Let $X_i$, $i=1,2$, be non-degenerate
continua such that each point of $X_i$ admits a closed
neighborhood with trivial $n_i$-dimensional cohomology, where
$n_i=\dim X_i$. If $X_1\times{}X_2$ is a quasi
$(n_1+n_2)$-manifold then each $X_i$ is a quasi $n_i$-manifold.

\endproclaim\ms

\noi{\bf Note}. Every polyhedron $P$ satisfies the condition from
this lemma: each point of $P$ admits a closed neighborhood which
is contractible.\qed\ms

\demo{Proof} Let $x_1\in X_1$ and $x_2\in X_2$ be arbitrary
points. We must construct open neighborhoods $V_i$ of $x_i$ in
$X_i$ satisfying the definition of a quasi $n_i$-manifold. Let
$n=n_1+n_2$. Since $X_1\times{}X_2$ is a quasi $n$-manifold, there
is an open neighborhood $V$ of the point $(x_1,x_2)$ in
$X_1\times{}X_2$ such that for every closed $(n{-}1)$-dimensional
set $F$ separating $X_1\times{}X_2$ between $(x_1,x_2)$ and
$(X_1\times{}X_2)\setminus{}V$ we have $H^{n-1}(F)\ne0$. Since the
same condition holds for any open neighborhood of $(x_1,x_2)$
contained in $V$, by Lemma 2D.2 and our hypothesis about $X_i$, we
may assume that $V=V_1\times{}V_2$, where each $V_i$ is an open
neighborhood of $x_i$ in $X_i$ such that
$H^{n_i}({\overline{V}_i})=0$, where ${\overline{V}_i}$ stands for
the closure of $V_i$ in $X_i$. We shall show that these $V_i$'s
are the desired neighborhoods.

By Lemma 2D.2 again, it follows that

\ms ($\ast$) {\it for any closed subset $A$ of $X_i$ contained in
${\overline{V}_i}$ and any $k\ge{}n_i$ we have $H^{k}(A)=0$}.\ms

\noi Now consider a closed $(n_i{-}1)$-dimensional subset $F_i$ of
$X_i$ separating $X_i$ between $x_i$ and $X_i\setminus{}V_i$. Then
$X_i\setminus F_i=U_i \cup W_i$, where $U_i, W_i$ are disjoint
open sets in $X_i$ such that $x_i\in U_i$ and $X_i\setminus
V_i\subset W_i$. Then $\overline{U}_i\subset{}V_i$, and the
boundary $\partial{}U_i$ ($\subset F_i$) separates $X_i$ between
$x_i$ and $X_i\setminus{}V_i$. It follows that
$\partial(U_1\times{}U_2)$ separates $X_1\times{}X_2$ between
$(x_1, x_2)$ and $(X_1\times{}X_2)\setminus{}V$. Thus
$H^{n-1}(\partial(U_1\times{}U_2))\ne0$.

Then consider the following portion of the Mayer-Vietoris
cohomology exact sequence of the couple $\{\partial{}U_1\times
\overline{U}_2,\overline{U}_1\times\partial{}U_2\}$:
$$H^{n-2}(\partial{}U_1\times\partial{}U_2)
@>{\delta^{\ast}}>> H^{n-1}(\partial(U_1\times{}U_2)) \to
H^{n-1}(\partial{}U_1\times \overline{U}_2)\oplus{}
H^{n-1}(\overline{U}_1\times\partial{}U_2)\ .$$ (The sequence
takes this form because $\partial{}U_1\times\partial{}U_2=
(\partial{}U_1\times
\overline{U}_2)\cap{}(\overline{U}_1\times\partial{}U_2)$ and
$\partial(U_1\times{}U_2)= (\partial{}U_1\times \overline{
U}_2)\cup{}(\overline{U}_1\times\partial{}U_2)$.) By the
K{\"u}nneth formula and ($\ast$), $H^{n-1}(\partial U_1 \times
\overline{U}_2)=0$ and
$H^{n-1}(\overline{U}_1\times\partial{}U_2)=0$. It follows that
$\delta^{\ast}$ is an epimorphism. Thus
$H^{n-2}(\partial{}U_1\times\partial{}U_2)\ne0$ since
$H^{n-1}(\partial(U_1\times{}U_2))\ne0$. Again, by the K{\"u}nneth
formula and ($\ast$), $H^{n-2}(\partial U_{1} \times
\partial{}U_2)= H^{n_1-1}(\partial{}U_1)\otimes
H^{n_2-1}(\partial{}U_2)$. It follows that both
$H^{n_1-1}(\partial{}U_1)$ and $H^{n_2-1}(\partial{}U_2)$ are not
trivial.

Note that $H^{n_i}(F_i,\partial{}U_i)=0$ since $\dim{}F_i=n_i-1$.
Thus, from the cohomology exact sequence of the pair
$(F_i,\partial U_i)$ and ($\ast$), it follows that the
homomorphism $H^{n_i-1}(F_i)\to{}H^{n_i-1}(\partial{}U_i)$ induced
by the inclusion $\partial{}U_i\hookrightarrow{}F_i$ is an
epimorphism. Consequently, each $H^{n_i-1}(F_i)$ is not trivial,
which concludes the proof of our lemma. \qed \enddemo

\proclaim{Lemma 2D.4} Let $p_i:X \to P_i$ and $q_i:S_i \to X$,
$i=1,\cdots,n$, be any mappings such that
$$H^1(p_i\circ q_j):H^1(P_i)\to H^1(S_j)$$ is an epimorphism
for each $i=j$, and the zero homomorphism for $i\ne j$. Then the
homomorphism
$$\varphi:H^1(X)\to H^1(S_1)\oplus\cdots\oplus{}H^1(S_n),$$ given by
$\varphi=(H^1(q_1),\cdots,H^1(q_n))$ is an epimorphism.
Consequently, $\rank{}H^1(X)\ge
\rank{}H^1(S_1)+\cdots+\rank{}H^1(S_n)$.
\endproclaim

\demo{Proof} Let $\beta=(\beta_1,\cdots,\beta_n)$ be any element
in $H^1(S_1)\oplus \cdots\oplus H^1(S_n)$. We have to find an
$\alpha\in H^1(X)$ such that $\varphi(\alpha)=\beta$. By our
assumption $\beta_i=H^1(p_i\circ q_i)(\alpha_i)$ for some
$\alpha_i\in H^1(P_i)$ and each $i$. One easily verifies that
$\alpha=H^1(p_1)(\alpha_1)+\cdots+H^1(p_n)(\alpha_n)$ satisfies
the equality. \qed \enddemo

In the following lemma we keep the notation and assumptions of
section 2C.

\proclaim{Lemma 2D.5} Let $M$ be a ramified $n$-manifold in the
product $P_1\times\cdots\times{}P_n$ of graphs. Then
$\rank{}H^1(M)\ge n$. If $\rank{}H^1(M)=n+k$, where $k<n$, we have
$n_{J_M}\geq n-k$. If $J_M=\{1,\cdots, n\}$ then $M$ coincides
with the $n$-torus $p_1(M)\times\cdots\times{}p_n(M)$; if $J_M$ is
a proper non-void subset of $\{1,\cdots, n\}$ then
$M=p_{J_M}(M)\times p_{J_M^c}(M)$, where $p_{J_M}(M)=\prod_{j\in
J_M}p_j(M)$ is an $n_{J_M}$-torus and $p_{J_M^c}(M)$ is a ramified
$n_{J_M^c}$-manifold which has no projection onto a circle.
\endproclaim

\demo{Proof} Let $v_j$, $j=1,\cdots,n$, denote a vertex of
$K'_{\{j\}^c}$. By Property (f), $P_j(v_j)$ is a finite union of
ramified 1-manifolds. To continue the proof we apply Lemma 2D.4 as
follows.

Let $q_j:P_j(v_j)\to{}M$ be the map such that $p_j\circ{}q_j$ is
the inclusion $P_j(v_j)\hookrightarrow{}P_j$ and
$(p_{\{j\}^c}\circ{}q_j)(P_j(v_j))=\{v_j\}$. Since $p_j\circ{}q_j$
is an inclusion into a graph, $H^1(p_j\circ{}q_j)$ is a
epimorphism. If $i\ne{}j$ then $H^1(p_i\circ{}q_j)$ is the zero
homomorphism as the image of $p_i\circ{}q_j$ is a point. Thus, by
Lemma 2D.4, we obtain

\ss($\ast$) $\rank{}H^1(M)\geq\rank{}H^1(P_1(v_1))+\cdots+\rank
H^1(P_n(v_n))$.

\ss \noindent Notice that $\rank{}H^1(P_j(v_j))\geq 1$. Hence
$\rank{}H^1(M)\geq n$, which proves the first assertion.

Now we prove the second one. To this end, pick the vertices $v_j$
so that $\rank{}H^1(P_j(v_j))$ $\geq \rank{}H^1(P_j(w_j))$ for
each vertex $w_j\in{}K'_{\{j\}^c}$. Let $$J_0= \{j\in
\{1,\cdots,n\}:\rank H^1(P_j(v_j))=1\}.$$ Since
$\rank{}H^1(M)=n+k$, by ($\ast$) we infer that
$\rank{}H^1(P_j(v_j))\geq2$ for at most $k$ indices $j$. Hence
$J_0$ consists of at least $n-k$ indices. If $j\in J_0$ then
$\rank{}H^1(P_j(w_j))=1$ for each vertex $w_j\in{}K'_{\{j\}^c}$.
Hence $P_j(w_j)$ is a circle. By Property~(h) we infer that
$p_j(M)$ is a circle for each $j\in{}J_0$. It follows that $J_0
\subset J_M$. Hence $n_{J_M}\geq n-k$. The remaining assertions
directly follow from Property (k). \qed\enddemo\ss

In section 2C we have established notation related to a subset $M$
lying in a product of graphs $P_1\times\cdots\times P_n$. Here we
introduce an analogous notation related to a subset $X$ of the
product of curves $Y_1\times\cdots\times Y_n$. In particular, for
a non-void subset $J$ of $\{1,\cdots,n\}$ let $Y_J$ denote the
product $\prod_{j\in J}Y_j$, and let $pr_J: Y_1\times\cdots\times
Y_n \to Y_J$ denote the corresponding projection. Similarly, let
$$J_X=\{j\in\{1,\cdots,n\}:pr_j(X)~ \text{\it {is a circle}}\}.$$
If $J_X$ is non-void then $X$ is said to {\it have projections
onto circles}, and any $pr_j$, for $j\in J_X$, is said to be {\it
a projection of X onto a circle}.

\proclaim{Theorem 2D.6} Let $X$ be a locally connected quasi
$n$-manifold in the product $Y_1\times\cdots\times{}Y_n$ of $n$
curves, $n\ge1$. Then we have: \ss{\rm(1)} $\rank{}H^1(X)\ge n$
$;$ \ss{\rm(2)} If $\rank{}H^1(X)=n+k$, where $k<n$, then $J_X$
contains at least $n-k$ elements. In particular, if
$\rank{}H^1(X)=n$ then $J_X=\{1,\cdots,n\}$ $;$ \ss{\rm(3)} If
$J_X=\{1,\cdots,n\}$ then $X=pr_1(X)\times\cdots\times pr_n(X)$ is
an $n$-torus $;$ \ss{\rm(4)} If $J_X$ is a proper non-void subset
of $\{1,\cdots,n\}$ then $X=(\prod_{j\in J_X}pr_j(X))\times
pr_{J_X^c}(X)$, where the first factor is an $n_{J_X}$-torus, and
$pr_{J_X^c}(X)$ is a quasi $n_{J_X^c}$-manifold in $Y_{J_X^c}$
having no projection onto a circle.
\endproclaim

\demo{Proof} Let
$f=(f_1,\cdots,f_n):X\to{}Y_1\times\cdots\times{}Y_n$ denote the
inclusion, i.e. $f_j=pr_j|X$ for each $j$. To prove (1) and (2) we
can assume that $\rank{}H^1(X)$ is finite. Then, by Theorem 2B.1,
there exist mappings $g
=(g_1,\cdots,g_n):X\to{}P_1\times\cdots\times{}P_n$ and
$h=h_1\times\cdots\times{}h_n:P_1\times\cdots\times{}P_n\to
Y_1\times \cdots\times Y_n$ such that $f=h\circ g$ (hence
$f_j=h_j\circ{}g_j$ for each $j$), where each $g_j:X\to P_j=|K_j|$
is a surjection onto a graph. Then $g$ is an embedding, hence
$M=g(X)$ is a quasi $n$-manifold in the product $P_1 \times \cdots
\times P_n$. So, $\rank{}H^1(X)=\rank{}H^1(M)$. By Theorem 2A.10,
$M=|K(M)|$, where $K(M)$ is a subcomplex of
$K_1\kw\cdots\kw{}K_n$. Then, by Proposition 2A.5(iii), we see
that $M$ is also a ramified $n$-manifold. Thus, by Lemma~2D.5,
$\rank{}H^1(X)\ge n$, which proves (1). By Lemma 2D.1(ii) we have
$J_X=J_M$. In fact, $h$ is a product of mappings $h_j:P_j\to Y_j$,
is injective on $M$ and $X=h(M)$, hence the hypotheses of that
lemma are fulfilled. Then the equality simply follows. Now we
shall prove (2). To this end, suppose $\rank H^1(X)=n+k$, where
$k<n$. Then $J_X(=J_M)$ contains at least $n-k$ elements, by Lemma
2D.5. This proves (2). Next we shall prove (3). So, suppose
$J_X=\{1,\cdots,n\}$. Then $J_M=\{1,\cdots,n\}$, hence by Lemma
2D.5 we have $M=p_1(M)\times\cdots\times p_n(M)$ (where
$p_j(M)=P_j$). Since $X=h(M)$, $h=h_1\times\cdots\times h_n$, and
$h_j(p_j(M)=pr_j(X)$ we infer that $X=(h_1\times\cdots\times
h_n)(p_1(M)\times\cdots\times p_n(M))=pr_1(X)\times\cdots\times
pr_n(X)$. This proves (3). Finally, we shall prove (4). So,
suppose $J_X$ is a proper non-void subset of $\{1,\cdots,n\}$.
Since $J_M=J_X$, by Lemma 2D.5 we have $M=p_{J_M}(M)\times
p_{J_M^c}(M)$, where $p_{J_M}(M)=\prod_{j\in J_M}p_j(M)$ with each
$p_j(M)$, $j\in J_M$, being a circle. By Lemma 2D.3 the set
$p_{J_M^c}(M)$ is a quasi $n_{J_M^c}$-manifold. Since $X=h(M)$,
$h=h_1\times\cdots\times h_n$, and $h|M$ is an embedding, this
factorization of $M$ and Lemma 2D.1 imply all assertions we need.
Indeed, Lemma 2D.1(i) implies that each $h_j$, $j\in J_M$, is an
embedding whose image equals $pr_j(X)$, and also
$h_{J_M^c}|p_{J_M^c}(M)$ is an embedding whose image equals
$pr_{J_X^c}(X)$. It follows that $X=(\prod_{j\in
J_X}pr_j(X))\times pr_{J_X^c}(X)$, where each $pr_j(X)$, $j\in
J_X$, is a circle. Since for no $j\in J_X^c$ the set $pr_j(X)$ is
a circle, the set $pr_{J_X^c}(X)$ has no projection onto a circle.
\qed\enddemo

\proclaim{Corollary 2D.7} Let $X$ be a locally connected quasi
$n$-manifold in the product $Y_1\times\cdots\times{}Y_n$ of $n$
curves, $n\ge1$. If no set $pr_j(X)$ is a circle then $\rank
H^1(X)\geq 2n$. \qed
\endproclaim

Applying the K\"{u}nneth formula we infer the following

\proclaim{Corollary 2D.8} No product $\T^{n-k} \times M^k$, where
$2\leq k\leq n$ and $M^k$ is a closed $k$-manifold with
$H^1(M^k)=0$, can be embedded in a product of $n$ curves. \qed
\endproclaim\ms

\noindent{\bf Remark.} This corollary implies the Borsuk theorem
\cite {B3}.\ms

The "Bing house" and the "dunce hat" are contractible quasi
$2$-manifolds. Hence the first cohomology group of both examples
is trivial. Thus the theorem implies that

\proclaim{Corollary 2D.9} Both the "Bing house" and the "dunce
hat" are $2$-dimensional compact contractible polyhedra, and
neither can be embedded in a product of two curves.\qed
\endproclaim\ms

\noindent {\bf Note.} This corollary shows that the number of
factors in the Nagata embedding theorem cannot be reduced to $n$,
even for contractible $n$-dimensional polyhedra.\qed\bs

\proclaim{Corollary 2D.10} Let $X$ be an n-manifold, $n\ge2$,
lying in the product $Y_1\times\cdots\times{}Y_n$ of $n$ curves.
If $\rank{}H^1(X)\le n+1$ then $X=S_1\times\cdots\times{}S_n$,
where each $S_j$ is a circle in $Y_j$.
\endproclaim

\noi{\bf Note.} The assumptions on $X$ can be relaxed: it is
enough to assume that $X$ is a locally connected quasi
$n$-manifold which is also a pseudo $n$-manifold.

\demo{Proof} By Theorem 2D.6 there is a set
$J\subset\{1,\cdots,n\}$ composed of $n-1$ elements such that

\ss $(\ast) \ \ X=(\prod_{j\in{}J}S_j)\times{}X'\ $,

\ss \noi where each $S_j$ is a circle in $Y_j$ and $X'$ is a quasi
$1$-manifold in $Y_i$, where $i$ is the element of the set
$\{1,\cdots,n\}\setminus{}J$. Since $X'$ is locally connected and
$\rank{}H^1(X')\le2$ it follows that $X'$ is a graph with no
endpoint. As $X$ is a pseudo $n$-manifold, by ($\ast$) it follows
that $X'$ contains no triod, hence it is a circle. This completes
the proof.
\enddemo

\centerline{2E. {\it Contractible $2$-dimensional polyhedra in
products of two graphs}}\ms

In this section we prove a result which in a particular case gives
a noteworthy property of 2-dimensional polyhedra acyclic in
dimension 1 and embeddable in products of two curves. As neither
the "Bing house" nor the "dunce hat" have this property, we get
another argument for non-embeddability of those examples in
products of two curves.

\proclaim{Theorem 2E.1} Let P be a $2$-dimensional connected
polyhedron embeddable in a product of two curves. If $P=|K|$,
where $K$ is a regular $CW$ complex, then

\roster \item "(i)" {\it if $H^1(P)=0$ then $K$ is collapsible};
\item "(ii)" {\it if $\rank H^1(P)=1$ then $K$ collapses to a
circle}; \item "(iii)" {\it if $\rank H^1(P)=2$ then $K$ collapses
to either a torus or a quasi $1$-manifold.\endroster
\endproclaim

\noindent{\bf Remark.} Also this theorem implies the Borsuk
theorem \cite {Bo3}.

\ms \demo{Proof} Let
$\kappa=\{K=K_0\searrow{}K_1\searrow\cdots\searrow{}K_n\}$ be a
maximal sequence of subcomplexes of $K$ such that each successive
complex is obtained from the preceding one by an elementary
collapsing. First we establish some general properties:

\ss (1) $K_n$ {\it is connected and $H^i(|K_n|)=H^i(P)$ for
$i\ge1$}. \ss

\ss (2) {\it If $K_n$ is $1$-dimensional then $|K_n|$ is a quasi
$1$-manifold, so $\rank H^1(|K_n|)\ge 1$}.\ss

\noi In fact, since $\kappa$ is maximal $K_n$ has no endpoint.

\ss (3) {\it If $K_n$ is $2$-dimensional and $\rank H^1(P)\le2$
then $|K_n|$ is a torus. So, $\rank H^1(|K_n|)=2$}.

\noi\ss In fact, let $X$ denote the union of all 2-cells of $K_n$.
Then $X=|K'_n|$, where $K'_n$ is a subcomplex of $K_n$. Since
$\kappa$ is maximal, each 1-cell of $K'_n$ is a face of at least
two different $2$-cells of $K'_n$. Thus each component of $X$ is a
ramified 2-manifold. Since $|K_n|\setminus{}X$ is 1-dimensional,
by exactness of the cohomology sequence of the pair $(|K_n|,X)$,
the homomorphism $H^1(|K_n|)\to H^1(X)$ induced by the inclusion
$X\hookrightarrow|K_n|$ is an epimorphism. It follows from the
assumption that $\rank{}H^1(X)\le2$. From Theorem~2D.6 we infer
that $X$ is homeomorphic to a torus. To complete the proof it
suffices to show that $X=|K_n|$. Suppose, to the contrary, that
$X$ is a proper subset of $|K_n|$. Let $C$ denote the closure (in
$|K_n|$) of a component of $|K_n|\setminus{}X$. Then $C$ is a
1-dimensional connected subpolyhedron of $|K_n|$ intersecting $X$
in a finite (nonzero) number of points. Since $\kappa$ is maximal,
each endpoint of $C$ belongs to $X$. Notice that if $C$ meets $X$
in one point then it contains a circle $S$, and if $C$ meets $X$
in at least two points then it contains an arc $L$ with end points
in $X$. In the first case $\rank{}H^1(X\cup{}S)=3$ and in the
second case $\rank{}H^1(X\cup{}L)=3$. It follows that
$\rank{}H^1(|K_n|)\ge3$, contrary to our assumption. This proves
the equality $X=|K_n|$, which completes the proof of (3).

The conclusion of our theorem readily follows from the above
properties. \qed
\enddemo

Theorem~2E.3 below is a partial converse of Theorem~2E.1. In the
proof of 2E.3 we need the following

\proclaim{Lemma 2E.2} Let $K_1$ and $K_2$ be regular
$1$-dimensional $CW$ complexes, and let $A$ be an oriented arc in
$|K_1|\times|K_2|$ which is a union of $1$-cells of $K_1\kw{}K_2$.
Then there exist regular $1$-dimensional $CW$ complexes
$K'_1\supset K_1$ and $K'_2\supset K_2$, and a disc $D \subset
|K'_1|\times|K'_2|$, such that

\roster

\item "(i)" each component of $|K'_i|\setminus|K_i|$ is a $1$-cell
of $K'_i$ with one endpoint removed,

\ss\item "(ii)" $D$ is a union of $2$-cells of $K'_1\kw{}K'_2$ and
$D\cap(|K_1|\times|K_2|)=A$.

\endroster
\endproclaim

\demo{Proof} Let $(v_0,w_0)$ denote the initial point of $A$.
Without loss of generality we may assume that $A$ can be presented
as a union of "vertical" and "horizontal" arcs as follows:
$$A=\{v_0\}\times w_0w_1\cup v_0v_1\times \{w_1\}\cup
\{v_1\}\times w_1w_2\cup v_1v_2\times \{w_2\}\cup\cdots$$ To
obtain $K'_1$ we enlarge $K_1$ adding mutually disjoint $1$-cells
$v_jv_j'$ standing out of $K_1$ (take $v_kv_k'=v_jv_j'$ if
$v_k=v_j$). Similarly, to obtain $K'_2$ we enlarge $K_2$ adding
mutually disjoint $1$-cells $w_jw_j'$ standing out of $K_2$. Hence
(i) holds. Consider the following discs (if they are defined):
\roster \itemitem{} $v_jv_j'\times{}w_jw_j',$

\ss\itemitem{}  $v_jv_j'\times{}w_jw_{j+1},$

\ss\itemitem{} $v_jv_j'\times{}w_{j+1}w_{j+1}',$ \ss\itemitem{}
$v_jv_{j+1}\times{}w_{j+1}w_{j+1}'$.
\endroster
\noindent If $A$ terminates at a point $(v_{n+1},w_{n+1})$,
$n\geq0$, then define $D$ to be the union of the discs for all
$j=0,\cdots,n$. If $A$ terminates at a point $(v_n,w_{n+1})$
define $D$ to be the union of the discs for all $j=0,\cdots,n-1$
and two initial discs for $j=n$. Then one can verify that $D$ is a
disc and condition (ii) holds as well.\qed
\enddemo

\proclaim{Theorem 2E.3} Let $K$ be a regular $2$-dimensional $CW$
complex. If $K$ is collapsible then $|K|$ is embeddable in a
product of two trees.
\endproclaim

\demo{Proof} We shall prove a stronger version of this theorem:

\ms(0){\it there exist trees $|K_1|$, $|K_2|$ and an embedding
$h:|K|\to|K_1|\times|K_2|$ such that $h(\sigma)$ is a union of
cells of $K_1\kw{}K_2$ for each cell $\sigma\in{}K$.}\ms

By our hypothesis there is a sequence of elementary collapses of
$K$ to a point $\star$:
$$K=L_n \searrow\cdots\searrow L_0=\{\star\}.$$ The proof of
(0) will be done ones we show that it holds for each $L_m$,
$m=0,\cdots, n$, in place of $K$. Obviously, if $m=0$ then (0) is
true. Now assume (0) holds for $m-1\geq 0$. It suffices to prove
it for $m$. By our assumption (0) holds for $L_{m-1}$. Hence there
exist an embedding $h':|L_{m-1}|\to|K_1|\times |K_2|$ as in (0).
Since $L_m \searrow L_{m-1}$ is an elementary collapsing, $|L_m|$
is a union of $|L_{m-1}|$ and $\tau$, where $\tau$ is either a
1-cell or a 2-cell of $L_m$. If $\tau$ is a 1-cell then
$|L_{m-1}|\cap \tau$ is a vertex $u_0\in L_{m-1}$. Then
$h'(u_0)=(v_1, v_2)$, where $v_i$ is a vertex of $K_i$. If $\tau$
is a 2-cell then $|L_{m-1}|\cap\tau$ is an arc $A'$ which is a
union of 1-cells of $L_{m-1}$, so the arc $A=h'(A')$ is a union of
1-cells of $K_1\kw{}K_2$. One easily sees that in the first case
$|L_m|$ embeds in $|K_1|\times |K'_2|$ as in (0), where $K'_2$ is
obtained from $K_2$ by adding a 1-cell $v_2v_2'$ standing out of
$|K_2|$. In the other case (0) follows from Lemma~2E.2. This ends
the proof. \qed
\enddemo

The final results of this section are devoted to embeddability of
cones over polyhedra into products.

\proclaim{Theorem 2E.4} Let $P$ be a $(k+l+1)$-dimensional
polyhedron, where $k,l\ge0$. Then there exist polyhedra $P'$ and
$P''$ with $\dim{}P'=k$ and $\dim{}P''=l$ such that the cone over
$P$ can be embedded in the product of cones over $P'$ and $P''$.
\endproclaim

This theorem is a consequence of the following two lemmas.

\proclaim{Lemma 2E.5} Let $P$, $k$ and $l$ be as in {\rm 2E.4}.
Then there exist polyhedra $P'$ and $P''$ with $\dim{}P'=k$ and
$\dim{}P''=l$ such that $P$ can be embedded in the join $P'*P''$.
\endproclaim

\demo{Proof} Let $P=|K|$, where $K$ is a simplicial complex. Put
$P'=|K^{(k)}|$, where $K^{(k)}$ is the $k$-skeleton of $K$. Define
$P''$ to be {\it the dual} to $P'$ in $P$, i.e. $P''$ is the union
of all simplices of the barycentric subdivision of $K$ which are
disjoint with $P'$. Then $\dim{P'}=k$ and $\dim{P''}=l$. Observe
that $P$ is pl isomorphic to a subpolyhedron of the join
$P'\ast{P''}$. This follows from the fact that for any simplex
$\sigma\in K$ with $\dim\sigma\ge{k}$ we have
$\sigma=\sigma'\ast\sigma''$, where $\sigma'$ is the $k$-skeleton
of $\sigma$ (with respect to the standard simplicial structure on
$\sigma$) and $\sigma''$ is the dual to $\sigma'$ in $\sigma$.
\qed
\enddemo

\proclaim{Lemma 2E.6} Let $P'$ and $P''$ be polyhedra. Then the
product of cones over $P'$ and $P''$ is homeomorphic to the cone
over the join $P'*P''$.
\endproclaim

\demo{Proof} For a polyhedron $Q$ let $aQ$ denote the cone with
vertex $a$ and base $Q$. According to the definition of link (cf.
\cite{R-S, p. 2}), $Q$ may be considered as a link of the vertex
$a$ in $aQ$. The conclusion of our lemma is the following special
case of a known formula (see \cite{R-S, p. 24, Ex. (3)}) for the
link of a vertex in a product of two polyhedra: $\lk((a', a''),
a'P'\times a''P'')\approx P'*P''$. So, $a'P'\times a''P'')\approx
c(P'*P'')$.\qed
\enddemo

\proclaim{Corollary 2E.7} The cone over an $n$-dimensional
polyhedron can be embedded in a product of $n+1$ copies of an
$m$-od.
\endproclaim

\demo{Proof} We prove this result by induction on the dimension
$n$. We start the induction with $0$-dimensional polyhedra, in
which case the proof is obvious. The inductive step is proven
applying Theorem 2E.4 for $k=0$ and $l=n-1$.\qed
\enddemo

Note that neither Lemma~2E.5 nor Corollary~2E.7 (and thus
Theorem~2E.4) can be extended to more general case of continua in
place of polyhedra. For example the Menger curve can not be
embedded in the join of two 0-dimensional compacta. Neither the
cone over the Menger curve is embeddable in the product of two
cones over 0-dimensional compacta.\ms

\centerline{2F. {\it On embeddings of $2$-dimensional polyhedra
into products of two graphes}} \centerline{-- {\it a solution of
Cauty's problem}}\bs

The main result of this section is Theorem 2F.6 which asserts that
there exist a 2-dimensional polyhedron which can be embedded in a
product of two curves but is not embeddable in any product of two
graphs. This provides a negative solution to a problem of R. Cauty
\cite{C}. The construction strongly depends on the following
special case of Theorem 2D.6: {\it Any topological torus in a
product of two curves is a product of two circles lying in
corresponding coordinate curves}. First we shall establish several
auxiliary lemmas.\ss

\proclaim{Lemma 2F.1} Let $Z_0, Z_1$ be two non-degenerate
compacta such that $Z_0\times Z_1$ is a topological cylinder. Then
one of $Z_i$ is a circle and the other is an arc.\qed
\endproclaim

Let us recall that by a $\theta$-{\it curve} we mean a union of
three arcs having common endpoints and else mutually disjoint. The
arcs are called {\it edges} and the endpoints - {\it vertices} of
the curve. Any $\theta$-curve can be regarded as the space of a
regular $CW$ complex with edges as 1-cells and vertices as
0-cells. This complex will be called {\it associated} with the
$\theta$-curve. Any $\theta$-curve contains three different
circles, each being a union of two different edges.

For any subsets $A_1,\cdots, A_n$ of a space $X$, and
$B_1,\cdots,B_n$ of a space $Y$ we have $$(A_1\times
B_1)\cap\cdots\cap (A_n\times B_n)=(A_1\cap\cdots\cap A_n)\times
(B_1\cap\cdots\cap B_n).$$ This formula is one of a few elementary
tools used in the proofs which follow.\ms

\centerline{\it From this point on, $P$ denotes a fixed
$\theta$-curve with vertices $a_0,a_1$,} \centerline {\it $K$
denotes the $CW$ complex associated with $P$, and} \centerline{\it
$S_i$, $i=0,1,2$, denote the circles in P.}\bs

 \proclaim{Lemma 2F.2} Let $h:P\times \Sb^1 \to Y_1\times Y_2$
be an embedding. Suppose $h(S_1 \times \Sb^1)=Q_1\times R_1$ and
$h(S_2\times \Sb^1)=Q_2\times R_2$, where $Q_1,Q_2$ are circles in
$Y_1$ and $R_1,R_2$ are circles in $Y_2$. Then either $(1)$
$Q_1\cup Q_2$ is a $\theta$-curve and $R_1=R_2$, or $(2)$
$Q_1=Q_2$ and $R_1\cup R_2$ is a $\theta$-curve.
\endproclaim

\demo{Proof} The intersection $(S_1\times \Sb^1)\cap (S_2\times
\Sb^1)=(S_1\cap S_2)\times \Sb^1$ is a topological cylinder. Hence
its image under $h$ is a topological cylinder as well. But
$h((S_1\times \Sb^1)\cap (S_2\times \Sb^1))=h(S_1\times \Sb^1)\cap
h(S_2\times \Sb^1)=(Q_1\cap Q_2)\times (R_1\cap R_2)$. Hence the
conclusion readily  follows from Lemma 2F.1. \qed\enddemo\ms

\proclaim{Lemma 2F.3} Let $h:P\times \Sb^1 \to P\times \Sb^1$ be a
homeomorphism. Then for every cell $\sigma\in K$ there exists a
cell $\sigma'\in K$ such that
$h(\kn\sigma\times\Sb^1)=\kn\sigma'\times\Sb^1$. Consequently,

$(1)$ $h(\sigma\times\Sb^1)=\sigma'\times\Sb^1$;

$(2)$ {\it for any $z\in \Sb^1$ the projection $pr_P$ maps
$h(\sigma\times\{z\})$ onto} $\sigma'$.
\endproclaim

\demo{Proof} Since $h$ maps the 2-manifold part
$M_2(P\times\Sb^1)$ onto itself, it permutes the components of
this set. Each component has the form $\kn\sigma\times\Sb^1$,
where $\sigma$ is a 1-cell of $K$. Hence the assertion follows for
1-cells. Likewise, $h$ maps $P\times\Sb^1\setminus
M_2(P\times\Sb^1)=\{a_0,a_1\}\times\Sb^1$ onto itself, hence it
permutes the components $\{a_0\}\times\Sb^1$ and
$\{a_1\}\times\Sb^1$. Hence the assertion holds for 0-cells as
well. Then (1) holds because $\cl
(\kn\sigma\times\Sb^1)=\sigma\times\Sb^1$, and (2) holds if
$\sigma$ is a 1-cell because $\sigma\times\{z\}$ connects
$\{a_0\}\times\Sb^1$ and $\{a_1\}\times\Sb^1$ (for
$\sigma\times\{z\}=\{(a_0,z)\} \;\cup \;\kn\sigma \times \{z\}
\;\cup \;\{(a_1,z)\}$), hence its image $h(\sigma\times\{z\})$
connects this circles as well. For 0-cells (2) follows from (1).
\qed
\enddemo\ms

\centerline{\it From this point on,} \centerline {\it $h:P\times
P\to Y_1\times Y_2$ denotes an embedding into a product of two
curves} \centerline{\it with $h(S_1 \times S_1)=Q_1\times R_1$,
$h(S_2 \times S_1)=Q_2\times R_1$,} \centerline{\it where
$Q_1,Q_2$ are circles in $Y_1$ and $R_1$ is a circle in $Y_2$.}
\bs

\proclaim{Lemma 2F.4} There is a circle $R_2$ in $Y_2$ such that
$h(S_1 \times S_2)=Q_1\times R_2$, $h(S_2 \times S_2)=Q_2\times
R_2$. Moreover, $Q_1\cup Q_2$ and $R_1\cup R_2$ are
$\theta$-curves.
\endproclaim

\demo{Proof} It follows from Lemma 2F.2 that $Q_1\cup Q_2$ is a
$\theta$-curve. Applying Lemma 2F.3 to the homeomorphism $P\times
S_1\to (Q_1\cup Q_2)\times R_1$ induced by $h$ we infer that
$pr_{Y_1}(h(S_i\times\{a_0\}))=Q_i$. So $pr_{Y_1}(h(S_i\times
S_2))\supset Q_i$ since $S_i\times\{a_0\}\subset S_i \times S_2$.
It follows that $h(S_1\times S_2)=Q_1\times R_2$ and $h(S_2\times
S_2)=Q_2\times R'_2$, where $R_2$ and $R'_2$ are circles in $Y_2$.
It follows from Lemma 2F.2 that $R_2=R'_2$. Now, apply Lemma 2F.2
to the embedding $S_1\times P\to Y_1\times Y_2$ induced by $h$.
Since $h(S_1\times S_1)=Q_1\times R_1$ and $h(S_1\times
S_2)=Q_1\times R_2$, we infer that $R_1\cup R_2$ is a
$\theta$-curve. This ends the proof. \qed\enddemo

Keeping the above assumptions, put $P_1=Q_1\cup Q_2$ and
$P_2=R_1\cup R_2$. By Lemma 2F.4 these are $\theta$-curves and
$h(P\times P)=P_1\times P_2$. Let $\{b_0,b_1\},\; \{c_0,c_1\}$
denote the vertices of $P_1,P_2$, respectively. Let $K_i$,
$i=1,2$, denote the regular $CW$ complex associated with $P_i$.

\proclaim{Lemma 2F.5} There exist homeomorphisms $h_i:P\to P_i$
such that h transforms each cell $\sigma_1\times\sigma_2\in K\kw
K$ onto $h_1(\sigma_1)\times h_2(\sigma_2)\in K_1\kw K_2$, i.e.
$h(\sigma_1\times\sigma_2)=h_1(\sigma_1)\times h_2(\sigma_2)$.
\endproclaim

\demo{Proof} We define $h_1$ by the formula:
$h_1(x)=pr_{Y_1}(h(x,a_0))$ for each $x\in P$. The mapping $h_2$
is defined similarly: $h_2(y)=pr_{Y_2}(h(a_0,y))$ for each $y\in
P$. It remains to show that these mappings have the desired
properties.

They are well defined because $h$ maps $P\times P$ into $P_1\times
P_2$. Next we shall show that $h_i$ maps $P$ homeomorphically onto
$P_i$. First we verify this for $h_1$. Applying Lemma 2F.3 to the
homeomorphism $P\times S_i\to P_1\times R_i$ induced by $h$, for
$i=1,2$, we infer  that for each $\sigma\in K$ there is
$\sigma'\in K_1$ such that $h(\kn\sigma\times
S_i)=\kn\sigma'\times R_i$. The assignment $\sigma\to\sigma'$ is
independent of $i$ because for any $\sigma\in K$ there is only one
$\sigma'\in K_1$ such that $h(\kn\sigma\times
\{a_0\})\subset\kn\sigma'\times R_i$. Since $P=S_1\cup S_2$ we
get\ss

(1) $h(\sigma\times P)=\sigma'\times P_2$.\ss

\noi Symmetrically, for each $\sigma\in K$ there is $\sigma''\in
K_2$ such that $h(S_i\times\kn\sigma)=Q_i\times\kn\sigma''$.
Hence, as above, we get\ss

(2) $h(P\times\sigma)=P_1\times\sigma''$.\ss

\noi In particular, $h(P\times\{a_0\})=P_1\times\{c_j\}$ for some
$j=0,1$. Thus, from the description of $h_1$ it follows that it is
a homeomorphism. Moreover, $h(\kn\sigma\times
\{a_0\})\subset\kn\sigma'\times P_2$. Therefore,
$h_1(\kn\sigma)\subset\kn\sigma'$ which implies
$h_1(\sigma)=\sigma'$. Likewise, $h_2$ is a homeomorphism and
$h_2(\sigma)=\sigma''$. It follows that
$h(P\times\sigma)=P_1\times h_2(\sigma)$ and $h(\sigma\times P)=
h_1(\sigma)\times P_2$.

Finally, combining (1) and (2), we infer that for any cells
$\sigma_1,\sigma_1\in K$ we have
$h(\sigma_1\times\sigma_2)=h((\sigma_1\times P) \cap
(P\times\sigma_2))=h_1(\sigma_1)\times h_2(\sigma_2)$, which
completes the proof.\qed\enddemo

\proclaim{Theorem 2F.6} Let $X=(P\times P)\cup D$, where $D$ is a
disc and $A=(P\times P)\cap D$ is an arc. Suppose the following
conditions are fulfilled:

{\rm(i)} $A\subset \partial D$$;$

{\rm(ii)} $(a_0,a_0)$ is an endpoint of A$;$

{\rm(iii)} $A\setminus\{(a_0,a_0)\}$ lies in the interior of a
$2$-cell of $K\kw K$.

\noi Then X is not embeddable in any product of two graphs though
it can be embedded in a product of two curves.
\endproclaim

\demo{Proof} In order to prove the first assertion suppose there
is an embedding $h:X\to Y_1\times Y_2$ into a product of two
graphs. Then, without loss of generality, we can assume that
$h(S_1 \times S_1)=Q_1\times R_1$ and $h(S_2 \times S_1)=Q_2\times
R_1$, where $Q_1,Q_2$ are circles in $Y_1$ and $R_1$ is a circle
in $Y_2$. (In fact, otherwise we obtain this assumption by
interchanging the position of coordinate spaces, see Lemma 2F.2.)
By Lemma 2F.5 there exist embeddings $h_i:P\to Y_i$ such that
$h(\sigma_1\times\sigma_2)=h_1(\sigma_1)\times h_2(\sigma_2)$ for
each cell $\sigma_1\times\sigma_2\in K\kw K$. Suppose
$A\setminus\{(a_0,a_0)\}$ lies in the interior of the $2$-cell
$\sigma\times\sigma\in K\kw K$. Since $Y_j$, $j=1,2$, is a graph
there is a neighborhood $U_j$ of $h_j(a_0)$ in $h_j(\sigma)$ such
that $U_j\setminus\partial h_j(\sigma)$ is an open subset of
$Y_j$. Then $U_1\times U_2$ is a neighborhood of $h(a_0,a_0)$ in
$h(\sigma\times\sigma)$. Hence there is a point $x\in
A\setminus\{(a_0,a_0)\}$ such that $h(x)\in(U_1\times U_2)\cap
h(A\setminus\{(a_0,a_0)\})$. Since
$x\notin\partial(\sigma\times\sigma)$ we get $h(x)\in (U_1\times
U_2)\setminus h(\partial(\sigma\times\sigma))=
(U_1\setminus\partial h_1(\sigma))\times (U_2\setminus\partial
h_2(\sigma))$. Since $x\in \cl(D\setminus(P\times P))$ we infer
that $h(x)\in \cl(h(D\setminus(P\times P)))$. Since $h(x)\in
(U_1\setminus\partial h_1(\sigma))\times (U_2\setminus\partial
h_2(\sigma))$ and this set is open in $Y_1\times Y_2$ there is a
point $y\in D\setminus(P\times P)$ such that $h(y)\in
(U_1\setminus\partial h_1(\sigma))\times (U_2\setminus\partial
h_2(\sigma))$. Hence $h(y)\in h(P\times P)$, a contradiction.

Now we shall show that $X$ can be embedded in the product $Y\times
Y$, where $Y$ is a curve (actually, $Y$ will be of order 3 and,
but one point, a local tree, see \cite{Kur} for definitions). To
construct $Y$ we proceed as follows. First of all, we fix a 1-cell
$\sigma\in K$ and pick a monotone sequence $a_2,a_3,\cdots$ of
different points lying in the interior of $\sigma$ and converging
to $a_0$. Define $Y$ to be the union $Y=P\cup a_2b_2\cup
a_3b_3\cup\cdots$, where $a_2b_2,a_3b_3,\cdots$ is a null-sequence
of mutually disjoint arcs (converging to $a_0$) such that each
$a_nb_n$ meets $P$ only at $a_n$. Next we define an arc $A^\star$
in $\sigma\times\sigma$ and a disc $D^\star$ in $Y\times Y$ by the
formulas:$$A^\star=\{(a_0,a_0)\}\cup\bigcup_{n\geq2}(\{a_n\}\times
a_na_{n+1}\cup a_na_{n+1}\times \{a_{n+1}\}),$$
$$D^\star=\{(a_0,a_0)\}\cup \bigcup_{n\geq2}[a_nb_n\times (a_nb_n\cup a_na_{n+1}\cup
a_{n+1}b_{n+1})\cup a_na_{n+1}\times a_{n+1}b_{n+1}].$$ Then
$D^\star\cap (P\times P)=A^\star$. Now we are ready to describe
the embedding. We do this in three steps.

First, define an embedding $(P\times P)\setminus
(\kn\sigma\times\kn\sigma)\to Y\times Y$ to be the inclusion.
Then, using elementary results, one can extend this embedding to
an embedding $P\times P\to Y \times Y$ which maps $\sigma \times
\sigma$ onto itself and $A$ onto $A^\star$. The resulting
embedding can be further extended to an embedding $X\to Y\times Y$
which maps $D$ onto $D^\star$. This completes the
proof.\qed\enddemo

\bs\centerline {\bf Problems}

\proclaim {Problem 2A.1} Is it possible to characterize a
polyhedron $|K|$ which is a quasi $n$-manifold in terms of the
complex $K$ itself?
\endproclaim

\proclaim {Problem 2B.1} Let $X$ be a locally connected quasi
$n$-manifold, $n\geq2$, with $H^1(X)$ of finite rank, and let
$(f_1,\cdots,f_n):X \to Y_1\times\cdots\times Y_n$ be an embedding
in a product of curves. Is it possible to approximate mappings
$f_i$ by mappings $f_i':X \to Y_i$ so that $(f_1',\cdots,f_n')$ is
still an embedding and each $f_i'(X)$ is a graph?\endproclaim

\proclaim{Problem 2D.1} Characterize quasi $n$-manifolds
embeddable in products of $n$ graphs.\endproclaim

\proclaim {Problem 2D.2}Let X be a locally connected pseudo
n-manifold lying in a product of n curves. Is X a locally
connected quasi n-manifold?
\endproclaim

\proclaim{Problem 2D.3} Characterize ordinary closed
$3$-manifolds. Must such a manifold be a product of non-degenerate
factors?\endproclaim

The following problem is of great interest.

\proclaim{Problem 2D.4} Let $X$ be a locally connected quasi
$n$-manifold lying in a product of $n$ curves. Does X admit an
essential map $X\to \Sb^n$?\endproclaim

\proclaim{Problem 2F.1} Characterize those $2$-dimensional
polyhedra which can be embedded into products of two graphs.

\endproclaim 

\enddocument